\RequirePackage{fix-cm}
\documentclass[smallextended]{svjour3}       
\smartqed  
\usepackage{graphicx}
\usepackage{amsmath,amssymb}
\usepackage{esint,amssymb}
\usepackage{graphicx}
\usepackage{mathtools}
\usepackage[colorlinks=true, pdfstartview=FitV, linkcolor=blue, citecolor=blue, urlcolor=blue,pagebackref=false]{hyperref}
\usepackage{microtype}
\usepackage{MnSymbol}
\newtheorem{Numproposition}{Numerical proposition}

\newcommand{\N}{\mathbb{N}}

\newcommand{\R}{\mathbb{R}}
\newcommand{\C}{\mathbb{C}}

\newcommand{\E}{\mathbb{E}}

\renewcommand{\tilde}{\widetilde}

\renewcommand{\hat}{\widehat}

\makeatletter
\newsavebox\myboxA
\newsavebox\myboxB
\newlength\mylenA
\newcommand*\mybar[2][0.75]{%
    \sbox{\myboxA}{$\m@th#2$}%
    \setbox\myboxB\null
    \ht\myboxB=\ht\myboxA%
    \dp\myboxB=\dp\myboxA%
    \wd\myboxB=#1\wd\myboxA
    \sbox\myboxB{$\m@th\overline{\copy\myboxB}$}
    \setlength\mylenA{\the\wd\myboxA}
    \addtolength\mylenA{-\the\wd\myboxB}%
    \ifdim\wd\myboxB<\wd\myboxA%
       \rlap{\hskip 0.5\mylenA\usebox\myboxB}{\usebox\myboxA}%
    \else
        \hskip -0.5\mylenA\rlap{\usebox\myboxA}{\hskip 0.5\mylenA\usebox\myboxB}%
    \fi}
\makeatother

\def\build#1_#2^#3{\mathrel{\mathop{\kern 0pt#1}\limits_{#2}^{#3}}}

\begin{document}

\title{Dynamical Fractional and Multifractal Fields
}


\author{Gabriel B. Apolin\'{a}rio        \and
        Laurent Chevillard   \and
Jean-Christophe Mourrat
}

\authorrunning{G. B. Apolin\'{a}rio, L. Chevillard and J.-C. Mourrat} 

\institute{G. B. Apolin\'{a}rio, L. Chevillard \at
              Univ Lyon, ENS de Lyon, Univ Claude Bernard, CNRS, Laboratoire de Physique, 46 all\'ee d'Italie F-69342 Lyon, France
              \email{gabriel.brito\_apolinario@ens-lyon.fr, laurent.chevillard@ens-lyon.fr}           
           \and
           J.-C. Mourrat \at
              Courant Institute of Mathematical Sciences, New York University, New York, New York, USA, and CNRS, France \email{jcm777@nyu.edu}
}

\date{Received: date / Accepted: date}

\maketitle

\begin{abstract}
Motivated by the modeling of three-dimensional fluid turbulence, we define and study a class of stochastic partial differential equations (SPDEs) that are randomly stirred by a spatially smooth and uncorrelated in time forcing term. To reproduce the fractional, and more specifically multifractal, regularity nature of fully developed turbulence, these dynamical evolutions incorporate an homogenous pseudo-differential  linear operator of degree 0 that takes care of transferring energy that is injected at large scales in the system, towards smaller scales according to a cascading mechanism. In the simplest situation which concerns the development of fractional regularity in a linear and Gaussian framework, we derive explicit predictions for the statistical behaviors of the solution at finite and infinite time. Doing so, we realize a cascading transfer of energy using linear, although non local, interactions. These evolutions can be seen as a stochastic version of recently proposed systems of forced waves intended to model the regime of weak wave turbulence in stratified and rotational flows. To include multifractal, i.e. intermittent, corrections to this picture, we get some inspiration from the Gaussian multiplicative chaos, which is known to be multifractal, to motivate the introduction of an
additional quadratic interaction in these dynamical evolutions. Because the theoretical analysis of the obtained class of nonlinear SPDEs is much more demanding, we perform numerical simulations and observe the non-Gaussian and in particular skewed nature of their solution.
\end{abstract}

\section{Introduction}
\label{s.intro}

The present investigation is mainly motivated by the modeling of some aspects of the random nature of fluid turbulence  \cite{TenLum72,Fri95}. To be more precise, let us begin with illustrating Kolmogorov's phenomenological theory of three-dimensional turbulence \cite{Kol41}. To do so, consider a component $u_{i\in\{1,2,3\}}(t,x)$, and $x\in\R^3$, of the divergence-free vector velocity field $u$ of a fluid of viscosity $\nu$, whose dynamics is governed by the incompressible three-dimensional Navier-Stokes equations. Moreover, we assume that this evolution is supplemented by an additive random vector forcing term $f(t,x)$ that we take divergence-free and smooth in space. Typically, without loss of generality, and to fix ideas, consider a zero-average, white-in-time Gaussian vector field whose covariance is of the form $\E[f(t,x)\cdot f(t',x')]\propto \delta(t-t')\mathcal C_f(|x-x'|)$, where the scalar  positive-definite function $\mathcal C_f(x)$ is $\mathcal C^{\infty}$ and takes significant values only for $|x|\le L$.

Experimental and numerical observations suggest that this dynamics, that we recall to be stirred all along the way by a random force, converges at large time towards a statistically stationary state in which the velocity variance $\sigma^2$ is finite and remains so in the fully developed turbulent regime (that is for $\nu\to 0$). Understanding how the fluid organizes itself spatially and temporally  to damp in an efficient way the energy that is injected at a given large scale $L$ is at the core of the phenomenology of turbulence. Indeed, a cascading process of energy is taking place, transferring in some ways the energy injected at the large scales $L$ to smaller scales, such that the fluid develops events of large spatial gradients and acceleration, that are eventually smoothed out by viscosity. As a consequence, in the asymptotic limit of infinite Reynolds number, or equivalently in the limit $\nu\to 0$, velocity becomes rough, as it can be quantified by the variance of the velocity increment $\delta_{\ell} u_i(t, x) = u_i(t,x + \ell)-u_i(t,x)$ and its decrease towards 0 according to
\begin{align}\label{eq:S2Turbu}
\lim_{\nu\to 0}\E \left( \delta_\ell u_i\right)^2\build{\sim}_{\ell \to 0}^{} c_2 \left( \frac{\ell}{L}\right)^{2H},
\end{align}
where $c_2$ is a positive constant of the order of $2\sigma^2$ and $H$ is the local H\"{o}lder exponent. In a turbulent context, it is universally observed that $H\approx 1/3$, as predicted by dimensional arguments \cite{TenLum72,Fri95}, and Eq. \ref{eq:S2Turbu} says that, at this second order statistical level, velocity shares the same local regularity as a fractional Brownian motion of Hurst parameter $H$ \cite{ManVan68}. Moreover, as a more precise characterization of the observed non-Gaussian nature of the velocity field, higher-order structure functions, i.e. the moments of order $q\in\N$ of the increments, behave as
\begin{align}\label{eq:SqTurbu}
\lim_{\nu\to 0}\E \left( \delta_\ell u_i\right)^q\build{\sim}_{\ell \to 0}^{} c_q \left( \frac{\ell}{L}\right)^{\zeta_q},
\end{align}
with a spectrum of exponents $\zeta_q$ which is possibly a nonlinear function of the order $q$. The deviation from the linear behavior $\zeta_q=qH$, which can be obtained starting from Eq. \ref{eq:S2Turbu}  and furthermore assuming that $u_i$ is a Gaussian field, is a manifestation of the intrinsically non-Gaussian nature of the fluctuations, known as the intermittency phenomenon, properly defined in the language of the multifractal formalism (see for instance \cite{Fri95,CheCas12}, and references therein).

In this spirit, the article is devoted to the design of a dynamics governed by a partial differential equation, forced by such a random force $f$, whose structure is simpler than the three-dimensional Navier-Stokes equations. The aim of this dynamics is to reproduce the aforementioned statistical properties of turbulence in the statistically stationary regime, without the ambition to mimic all of the behaviors inherent to fluid mechanics, such as laminar flows and the mechanisms of transition towards turbulence. To further simplify this picture, we will limit ourselves to one-dimensional space $x\in\R$, and consider a unique velocity component $u(t,x)\in \C$, the imaginary nature of such a modeled velocity will become clear later when energy conservation is discussed.

Reproducing a cascading process of energy from the large to the small scales is the great success of shell models (see the review articles \cite{BohJen98,Bif03}, and \cite{ConLev06,CheFri10,BarMor11} for a more mathematically inclined approach). They consist in considering a coupled system of nonlinear ordinary differential equations, inspired by the expression of the Navier-Stokes equations in Fourier space, each of them governing the evolution of a shell $u_n(t)\in\C$ with $n\in\N$, which is meant to mimic some aspects of the behavior of a velocity Fourier mode $\hat{u}(t,k_n)$ over a logarithmically-spaced lattice $k_n=k_02^n$. The dynamics is furthermore supplemented by a viscous damping term and a forcing term $f_n$ that is restricted to large scales (i.e. $f_n=0$ for say $n\ge 3$). The coupling of each shell $u_n$ with its closest neighbors, typically at larger scales $(u_k)_{n-2\le k\le n}$ and at smaller scales  $(u_k)_{n\le k\le n+2}$, is made in a heuristic and nonlinear way such that for instance the dynamics preserves some key invariants that share the same structure as kinetic energy and helicity. Shell models can be viewed as a dynamical system that possesses as many degrees of freedom as the number of shells, once boundary conditions are set in an appropriate way. In a certain sense, such shells are observed numerically to behave in a similar way as in turbulence (Eqs. \ref{eq:S2Turbu} and \ref{eq:SqTurbu}) \cite{Bif03}, but a complete analytical understanding of the energy transfer mechanisms remains an open question \cite{ConLev06}. In this spirit, it is shown in Ref. \cite{MatSui07} that, instead of considering a nonlinear coupling between the shells, a peculiar linear coupling that mimics a derivative with respect to the number of the shell $n$ is able to reproduce some aspects of the cascading process of energy. We will later employ this idea, which can be viewed as a transport equation in the scale-space, that can be fully understood since only linear interactions are considered.

Although the underlying idea of shell models is appealing, it is not clear how to interpret such a shell $u_n(t)$. Indeed, it is has been observed in various direct numerical simulations of the Navier-Stokes equations \cite{BruPum01} and in experiments \cite{CheMaz05} that the real and imaginary parts of the true Fourier modes are mostly Gaussian, being compatible with Eq. \ref{eq:S2Turbu} but not with  Eq. \ref{eq:SqTurbu}. Although this observation makes shells not clearly related to Fourier modes over a logarithmically-spaced lattice, some ways to interpret shells in a continuous framework are proposed in Ref. \cite{Mai15}, giving a meaning to the shells as Fourier modes, allowing to design related partial differential equations in physical space.

As we can see, on the one hand, interpreting shells as Fourier modes is not fully satisfactory. On the other hand, it is tempting to interpret shells as wavelet coefficients in a dyadic decomposition of velocity over a tree. Such a decomposition can been shown to be orthonormal for square-integrable functions, and possesses a reconstruction formula in physical space \cite{Dau92}, that has been used in a turbulent context \cite{BenBif93,ArnBac98a} in order to synthesize random fields able to reproduce aforementioned statistical properties (Eqs. \ref{eq:S2Turbu} and \ref{eq:SqTurbu}). Doing so, inverting this orthonormal decomposition in order to get the dynamics in the physical space requires to link these shells, or wavelet coefficients, both in scale and in space. This interpretation of shell models, much more complete than only considering interactions through scales, has been already explored in the literature \cite{BarBia13,BiaMor17}. Nevertheless, an analytical derivation of the statistical properties of such shells when the dynamics is forced by an external large-scale forcing remains difficult since further relations between these coefficients in space must be prescribed, in order
to guarantee, for instance, spatial homogeneity of the velocity field in physical space, i.e. that the underlying probability law is invariant by translation. Designing such an interaction between the shells is not obvious and barely discussed. Up to now, we are not aware of such a dyadic model over a tree able to reproduce the rough behaviors depicted by the behaviors of structure functions at small scales (Eqs. \ref{eq:S2Turbu} and \ref{eq:SqTurbu}) in a statistically stationary and homogeneous framework.
It is also worth mentioning the approach of the so-called Leith \cite{leith1967,thalabard2015} and EDQNM \cite{Ors70,BosChe12} models in directly proposing a PDE for the energy spectrum, based on the phenomenology of turbulence. Nevertheless, these models do not address fundamental statistical features of the underlying velocity field, such as homogeneity, stationarity and intermittency.

In a very different context, devoted to the mathematical understanding of the phenomenon of convergence of internal \cite{MaaLam95} and inertial \cite{RieVal97} waves towards attractors, as observed experimentally in linearly stratified flows \cite{MaaBen97,ScoErm13,BroErm16}, the authors of Ref. \cite{ColSai20} propose an original interpretation. They show that these waves, whose dispersion relation between their wavelength and their frequency is very peculiar, can be obtained as solutions of a linear partial differential equation (PDE), supplemented by an external forcing, where enters an homogeneous operator of degree 0 \cite{Col20,DyaZwo19}. Furthermore, and as a consequence of the nature of this operator, the phenomenon of attraction of waves is seen as a cascading process \cite{ColSai20}. As we will see, this operator can be interpreted as a linear transport in the Fourier space, and interestingly, its discretized version coincides with the linear shell model developed in  Ref. \cite{MatSui07}. It thus becomes very tempting to include such an operator in a dynamics that would transfer energy from the large to the small scales, as demanded by the phenomenology of turbulence. Doing so, this would mean that this cascading process could be captured by a linear mechanism. This is what we propose to study in the present article.

To go further in the presentation of our results, let us consider a one-dimensional velocity field $u(t,x)\in\C$ with $x\in \R$, and its continuous Fourier transform
\begin{align}\label{eq:FTu}
\mathcal F[u](t,k)\equiv \hat{u}(t,k)=\int_{x\in\R}e^{-2i\pi kx}u(t,x)dx.
\end{align}
In the sequel, we will be studying the following nonlinear stochastic partial differential equation given by
\begin{align}\label{eq:IntroPDE}
\partial_tu_{H,\gamma,\nu}&=P_H\mathcal LP_H^{-1}u_{H,\gamma,\nu}\\
&+ \gamma P_H\left[\left(\tilde{P}_0\mathcal L P_H^{-1}u_{H,\gamma,\nu}\right) \left( P_H^{-1}u_{H,\gamma,\nu}\right) \right] +\nu \partial^2_x u_{H,\gamma,\nu} + f,\notag
\end{align}
where we use the notation $\partial_t\equiv \partial/\partial t$ for temporal  and  $\partial_x\equiv \partial/\partial x$ for spatial derivatives. Viscosity $\nu$ enters in the dynamics through the second-order spatial derivatives $ \partial^2_x\equiv \partial^2/\partial x^2$. Henceforth, the forcing term $f(t,x)\in\C$ will be assumed Gaussian and uncorrelated in time, statistically homogeneous, with zero average and covariance given by
\begin{align}\label{eq:CovGaussForce}
\E\left[ f(t,x)f^*(t',x')\right]=\delta(t-t')\mathcal C_f(x-x'),
\end{align}
where $^*$ stands for the complex conjugate, and $\mathcal C_f$ is a smooth function that decays rapidly away from the origin. To fully determine the forcing $f\in\C$, we furthermore take $\E\left[ f(t,x)f(t',x')\right]=0$, which implies that its real and imaginary parts are chosen independently. For analytical and numerical purposes, we will for instance consider $\mathcal C_f(x)=\exp\left(-\frac{x^2}{2L^2}\right)$, where the large length scale $L$ will eventually coincide with the correlation length scale of velocity $u$, and is known in turbulence phenomenology as the integral length scale.

Several operators and parameters enter in the dynamics of the velocity field $u(t,x)$ (Eq. \ref{eq:IntroPDE}). Let us begin with the operator $\mathcal L$ that, as we will see, is responsible for the transfer of energy from the large scale $L$ towards smaller ones. The crucial step, made in Ref. \cite{MatSui07} in a discrete setup related to the dynamics of a shell model and in Ref. \cite{ColSai20} in a continuous one insightfully related to the propagation of waves in rotational and stratified flows, lies in demonstrating that such a transfer of energy through scales can be done in a linear fashion. In our continuous set up, we thus consider the linear operator
\begin{equation}\label{eq:IntroDefL}
\mathcal Lu(t,x) \equiv 2i\pi cxu(t,x),
\end{equation}
where $c$ is a constant. Using the language developed in Refs. \cite{ColSai20,Col20,DyaZwo19}, we can say that $\mathcal L$ is an homogeneous pseudo-differential operator of degree 0. From a physical point of view, the picture gets very clear in Fourier space, while writing Eq. \ref{eq:IntroDefL} in a equivalent way as
\begin{equation}\label{eq:IntroDefFourL}
\mathcal F[\mathcal Lu](t,k) =-c\partial_k\hat{u}(t,k),
\end{equation}
which says that the inviscid and unforced dynamics $\partial_tu=\mathcal Lu$ is nothing else than a transport equation in the Fourier space, i.e. $\partial_t\hat{u}=-c\partial_k\hat{u}$, towards increasing wavelengths for positive rate $c>0$, and respectively towards decreasing wavelengths for $c<0$. For the sake of clarity, let us consider $c>0$ such the transport goes in the direction of increasing $k$, as it is observed for turbulence. We will show in the sequel that once sustained by a forcing term $f$ (Eq. \ref{eq:CovGaussForce}), this intermediate dynamics will generate a solution $u_0(t,x)$, with initial condition $u_0(0,x)=0$, that eventually behaves similarly as a complex Gaussian white noise in space as $t\to \infty$, in a way that we examine during the course of the article.

In other words, the operator $\mathcal L$ entering in the full dynamics written in Eq. \ref{eq:IntroPDE} participates in transferring the energy injected at the scale $L$ to smaller scales, in a linear and non dissipative way, such that the solution $u_0(t,x)$ seen as a function of $x$ at a given large time $t$ develops the regularity of a white noise. Let us keep in mind that our aim for $u$ is to reproduce instead, at least from a statistical point of view (Eq. \ref{eq:S2Turbu}), the regularity of a fractional Gaussian field of parameter $H\in]0,1[$. For this purpose, we introduce the operator $P_H$ in the dynamical evolution which reads,
\begin{equation}\label{eq:IntroDefPH}
P_Hu(t,x) \equiv \int e^{2i\pi k x}\frac{1}{|k|_{1/L}^{H+1/2}}\hat{u}(t,k)dk,
\end{equation}
where a regularized absolute value $|\cdot |_{1/L}$ over the wavelength $1/L$ is introduced, such that  $|k |_{1/L}\approx |k|$ when $|k|\gg 1/L$ and $|k |_{1/L}\approx 1/L$ when $|k|\ll 1/L$. The inverse $P_H^{-1}$ of this operator reads accordingly
\begin{equation}\label{eq:IntroDefPHinverse}
P_H^{-1}u(t,x) \equiv \int e^{2i\pi k x}|k|_{1/L}^{H+1/2}\hat{u}(t,k)dk.
\end{equation}
We will see that the linear part of the full dynamics (Eq. \ref{eq:IntroPDE}), that is $\partial_tu_{H,0}=P_H\mathcal LP_H^{-1}u_{H,0}+ f$, eventually generates a solution $u_{H,0}(t,x)$, with initial condition $u_{H,0}(0,x)=0$, seen as a function of space $x$ and at a fixed and large time $t$, that shares several properties with a statistically homogeneous fractional Gaussian field, again as $t\to\infty$. In particular, the variance of $u_{H,0}$ will reach a finite value and the second order structure function will behave as in Eq. \ref{eq:S2Turbu}. We will also see how an additional viscous term generating a solution noted $u_{H,\nu}$ modifies this picture and  allows to reach furthermore a statistically stationary state.

Ultimately, let us discuss the nonlinear part of the dynamics  (Eq. \ref{eq:IntroPDE}) that we are proposing. Whereas all the ingredients that we previously discussed are based on linear operations on the velocity field $u$, and as we will see can be fully understood on a rigorous ground, this additional nonlinearity makes the overall dynamics much more intricate. For this reason, we will mostly rely on numerical simulations and observe their implications.

In a few words, the structure of this nonlinearity originates from the probabilistic construction of multifractal random fields using the Gaussian multiplicative chaos \cite{Man72,Kah85,RhoVar14}. In our setup, to be more precise, and as it is demanded by the complex nature of our dynamics (Eq. \ref{eq:IntroPDE}), we will invoke a complex generalization of this probabilistic object, some aspects of which have been already explored in the literature \cite{LacRho15}. It consists in introducing a random field able to reproduce key ingredients that enter in the nonlinear behavior of the spectrum of exponents $\zeta_q$ of high-order structure functions (Eq. \ref{eq:SqTurbu}). It is obtained as the exponential of a logarithmically correlated Gaussian random field. It can be seen as a particular case of the more general class of log-infinitely divisible measures \cite{BarMan02,SchMar01,BacMuz03,ChaRie05,RhoSoh14,RhoVar14} and has been extensively used under various forms while modeling the random nature of fluid turbulence \cite{SchLov87,BacDel01,MorDel02,RobVar08,CheRob10,PerGar16,CheGar19}. Moreover, the Gaussian field entering in the construction, that we recall to be logarithmically correlated, can be seen, in a way that we will discuss later, as a fractional Gaussian field of vanishing parameter $H=0$ \cite{DupRho17}. Not only is this remark important because such a field can thus be defined as a solution of regularized versions of random walks \cite{AloNua00,Che17,CheLag20}, but also because it makes a clear connection with the aforementioned build-up of fractional Gaussian fields using the operator $P_H$ (Eq. \ref{eq:IntroDefPH}) for the boundary case $H=0$. For several reasons that are developed in the sequel during the ad-hoc construction of multifractal fields, it turns out necessary to introduce a Hermitian symmetric version $\tilde{P}_0$ of the Fourier multiplier of such an operator $P_0$ (Eq. \ref{eq:IntroDefPH}) that  reads
\begin{equation}\label{eq:IntroDefP0tilde}
\tilde{P}_0u(t,x) \equiv -i \int e^{2i\pi k x}\frac{k}{|k|_{1/L}^{3/2}}\hat{u}(t,k)dk.
\end{equation}
Developing on these ideas, the design of the nonlinear term of (Eq. \ref{eq:IntroPDE}) is a consequence of rules of construction of multifractal fields able to reproduce the behaviors of the second (Eq. \ref{eq:S2Turbu}) and high-order (Eq. \ref{eq:SqTurbu}) structure functions. Interestingly, this approach which is based on a probabilistic ansatz also gives a way to define a multiplicative chaos as a solution of a dynamical process governed by a partial differential equation, forced by a smooth term, in a different spirit and setup than those proposed in the context of Liouville measures and two-dimensional quantum gravity \cite{Gar20,DubShe21}. We will see that, in our context, this probabilistic ansatz, not only for the multiplicative chaos, but more appropriately for a multifractal velocity field $u$, suggests a dynamics that is not closed in a simple fashion in terms of this field $u$. We then rely on a closure approach to simplify this dynamics, ending up with the quadratic nonlinearity that enters in Eq. \ref{eq:IntroPDE}.

The last parameter $\gamma\in\R$ entering in the proposed dynamics (Eq. \ref{eq:IntroPDE}) has the same origins. Its role in the aforementioned multifractal probabilistic ansatz is clear and governs entirely the level of multifractality and related non-Gaussian behaviors. Once inserted in our dynamics, a fully rigorous approach is much more demanding. Instead, we propose and design numerical simulations of the dynamics (Eq. \ref{eq:IntroPDE}) that indeed show that $\gamma$ governs, among others, the non-Gaussian nature of the solution $u(t,x)$, at least for the range of values that we have explored.

\smallskip
\subsection*{Organization of the paper.} We develop in Section \ref{Sec:HamiltonianL} the Hamiltonian dynamics induced by the homogenous pseudo-differental operator $\mathcal L$ (Eq. \ref{eq:IntroDefL}) of degree 0. More precisely, we compute the statistical properties at large time of the solution $u_{\nu}(t,x)$ of the partial differential equation $\partial_tu_{\nu}=\mathcal Lu_{\nu}+\nu\partial_x^2u_{\nu}+f$, with and without the additional action of viscosity $\nu$ and forcing $f$. We focus in Section \ref{Sec:LinearFractDyn} on the linear part of our proposed dynamics (Eq. \ref{eq:IntroPDE}), induced by the joint action of the operators $\mathcal L$ (Eq. \ref{eq:IntroDefL}) and $P_H$ (Eq. \ref{eq:IntroDefPH}). To do so, we first recall in paragraph \ref{Par:fGf} some key ingredients of the construction of fractional Gaussian field as linear operations on a Gaussian white noise measure, and then in paragraph \ref{Par:IndFracDynfGf} move on to the calculation of the statistical properties and regularity of the solution of the linear stochastic PDE $\partial_tu_{H,\nu}=P_H\mathcal LP_H^{-1}u_{H,\nu}+\nu\partial_x^2u_{H,\nu}+f$, again with and without the additional action of viscosity $\nu$ and forcing $f$. We then develop in Section \ref{Sec:NonLinearMultiDyn} the extension of the linear approach in order to go beyond the Gaussian framework. We begin in paragraph \ref{Par:ComplexMC} by recalling some basic facts about multifractal random fields, and develop a method for their construction starting from an instance of $u_0(t,x)$, which is similar to a complex Gaussian white noise at infinite time. Based on this method of construction, which is viewed as a probabilistic ansatz, we develop in paragraph \ref{Par:InducedDynComplexMC} the induced dynamics. Finally, because a rigorous approach aimed at calculating the statistical properties of the solution of the proposed stochastic PDE (Eq. \ref{eq:IntroPDE}) is for now out of reach, we design in Section \ref{Sec:Nums} a numerical algorithm and run simulations that give access to the solution $u_{H,\gamma,\nu}(t,x)$ in the statistical stationary regime. This allows us to estimate its statistical properties for a given set $(H,\gamma)$ of parameters, while considering averages across space of an instance at large time of the spatial profile $u_{H,\gamma,\nu}(t,x)$. As we will see, indeed the solution of Eq. \ref{eq:IntroPDE} reproduces the statistical properties announced in Eqs. \ref{eq:S2Turbu} and \ref{eq:SqTurbu}, and moreover exhibit an interesting non-vanishing third-order moment of the increments of the real part of $u_{H,\gamma,\nu}$.

\section{Hamiltonian dynamics induced by the cascade operator $\mathcal L$} \label{Sec:HamiltonianL}

The purpose of this section is the presentation of the first ingredient entering in the dynamics under study (Eq. \ref{eq:IntroPDE}), which concerns the statistical properties of the solution $u_\nu(t,x)\in \C$  of the following stochastic PDE
\begin{align}\label{eq:DynLnuF}
\partial_tu_\nu=\mathcal Lu_\nu+\nu\partial_x^2u_\nu+f,
\end{align}
where $\nu\ge 0$ is the viscosity, $f$ a Gaussian random forcing term, whose covariance is given in Eq. \ref{eq:CovGaussForce}, and the linear operator $\mathcal L$ (Eq. \ref{eq:IntroDefL}) that we recall the expression for any function $u$ for convenience,
\begin{equation}\label{eq:SecDefL}
\mathcal Lu(t,x) \equiv 2i\pi cxu(t,x).
\end{equation}
Without the forcing term $f$, observe that the evolution of the field $v(t,x)\equiv u(t,-x)$ is the same as in Eq. \ref{eq:DynLnuF} but with opposite rate $-c$ entering in the expression of $\mathcal L$ (Eq. \ref{eq:SecDefL}). Without viscous diffusion and forcing, the proposed dynamics (Eq. \ref{eq:DynLnuF}) is Hamiltonian, as can be seen from the skew-Hermitian symmetry of the operator $\mathcal L$ (Eq. \ref{eq:SecDefL}), and it preserves energy, i.e. the energy budget of the solution $u_0(t,x)$ is simple and given by $\partial_t|u_0|^2=0$.

\begin{proposition}\label{Prop:HamDynL}\textit{(Concerning the Hamiltonian dynamics induced by $\mathcal L$})

Consider the evolution
\begin{align}\label{eq:PDEHamilL}
\partial_t u_{0}(t,x)= \mathcal Lu_{0}(t,x)+ f(t,x),
\end{align}
where $f(t,x)$ is a Gaussian random force defined in Eq. \ref{eq:CovGaussForce}, with $\mathcal C_f$ a real and even function of its argument, and the linear operator $\mathcal L$ defined in Eq. \ref{eq:SecDefL}. Starting from the initial condition $u_0(0,x)=0$, the solution of this evolution is statistically homogeneous, meaning that the correlation in space
\begin{align}\label{eq:ExprCorrUXMY}
\mathcal C_{u_0}(t,x-y)= \mathbb E[u_{0}(t,x)u_{0}^*(t,y)],
\end{align}
is a function of the difference $x-y$. As a consequence of choosing independently the real and imaginary parts of the forcing $f$, we also have the following property at any time and positions,
\begin{align}\label{eq:ExprCorrUReIm}
\E[u_{0}(t,x)u_{0}(t,y)]=0.
\end{align}
Furthermore, $u_0(t,x)$ behaves similarly to a white noise at large time, such that, for any smooth function $g$,
\begin{align}\label{eq:BehavCuWN}
\lim_{t\to \infty}\int g(x)\mathcal C_{u_0}(t,x)dx= \frac{\mathcal C_f(0)}{2|c|} g(0).
\end{align}
In this sense, we would say that the operator $\mathcal L$ has transferred the energy from the large scale $L$ at which it is injected by the force $f$ towards vanishing scales at infinite time.
\end{proposition}

\begin{remark}
The proof of Proposition \ref{Prop:HamDynL} is straightforward and is a consequence of the exact expression of the solution of Eq. \ref{eq:PDEHamilL}. We have, starting with $u_0(0,x)=0$,
$$ u_0(t,x)=\int_0^t e^{2i\pi cx(t-s)}f(s,x)ds,$$
such that
$$\mathcal C_{u_0}(t,x-y)=\mathcal C_f(x-y) \frac{e^{2i\pi c(x-y)t}-1}{2i\pi c(x-y)},$$
where the function $\mathcal C_{u_0}(t,x)$ is defined in Eq. \ref{eq:ExprCorrUXMY}. To clarify its behavior at large time $t>0$, integrate it against a smooth function $g$ and obtain
\begin{align*}
\int g(x)\mathcal C_{u_0}(t,x)dx&=\int g(x)\mathcal C_f(x) \frac{e^{2i\pi cxt}-1}{2i\pi cx}dx\\
&=\frac{1}{|c|}\int g\left(\frac{x}{ct}\right)\mathcal C_f\left(\frac{x}{ct}\right) \frac{e^{2i\pi x}-1}{2i\pi x}dx\\
&\build{\sim}_{t\to\infty}^{}\frac{g(0)\mathcal C_f(0)}{|c|}\int  \frac{e^{2i\pi x}-1}{2i\pi x}dx,
\end{align*}
where the remaining indefinite integral is equal to $1/2$, which entails Eq. \ref{eq:BehavCuWN}.
\end{remark}

\begin{remark}
Taking into account a finite viscosity $\nu>0$ in this picture and then solving Eq. \ref{eq:DynLnuF} instead of Eq. \ref{eq:PDEHamilL} is also straightforward. We get for the same vanishing initial condition $u_\nu(0,x)=0$, using the Fourier transform defined in Eq. \ref{eq:FTu},
$$ \hat{u}_\nu(t,k)=\int_0^t e^{-(2\pi)^2\nu \left[ ks(k-cs)+\frac{c^2s^3}{3}\right]}\hat{f}(s,k-cs)ds.$$
Then, by taking expectations, we obtain an expression for the covariance function of $u_{\nu}$ (Eq. \ref{eq:ExprCorrUXMY}),
$$ \hat{\mathcal C}_{u_\nu}(t,k) =\int_0^t e^{-2(2\pi)^2\nu \left[ks(k-cs)+\frac{c^2s^3}{3}\right]} \hat{\mathcal C}_f(k-cs)ds,$$
and we recover the results obtained in the Hamiltonian case developed in Proposition \ref{Prop:HamDynL} while considering $\nu=0$.
Contrary to the inviscid case (i.e. $\nu=0$), the statistical behavior of the solution $u_{\nu}$ is very different when viscosity is finite. In particular, the variance reaches a finite value given by
\begin{align}
\lim_{t\to\infty}\E|u_\nu(t,x)|^2=\int_{s\in\R^+}\int_{k\in\R} e^{-2(2\pi)^2\nu \left[ks(k+cs)+\frac{c^2s^3}{3}\right]} \hat{\mathcal C}_f(k)dsdk.\notag
\end{align}
To see how it behaves as $\nu\to 0$, rescale the dummy variable $s$ by $\nu^{-1/3}$ and get
\begin{align}
\lim_{t\to\infty}\E|u_\nu(t,x)|^2&=\frac{1}{\nu^{1/3}}\int_{s\in\R^+}\int_{k\in\R} e^{-2(2\pi)^2\nu  \left[ks\nu^{-1/3}(k+cs\nu^{-1/3})+\frac{c^2s^3}{3\nu}\right]} \hat{\mathcal C}_f(k)dsdk,\notag\\
&\build{\sim}_{\nu\to 0}^{}\frac{\mathcal C_f(0)}{\nu^{1/3}|c|^{2/3}}\int_{s\in\R^+} e^{-\frac{2}{3}(2\pi)^2  s^3} ds,\label{eq:EquivLimVarUHamilt}
\end{align}
where the remaining integral can be evaluated with the help of the Gamma function. By inspection of the limiting behavior given in Eq. \ref{eq:EquivLimVarUHamilt}, we can see that the variance of $u$ is proportional to the one of the forcing term $\mathcal C_f(0)$, weighted by the diverging factor $\nu^{-1/3}$ as $\nu$ goes to zero.
\end{remark}

\begin{remark}
As we can see, in the presence of viscosity, the solution of Eq. \ref{eq:DynLnuF} reaches a statistically steady state, in which the variance is finite. Let us underline that the role of the transfer term played by $\mathcal L$ is crucial in the establishment of this regime. Indeed, without this term, i.e. taking $c=0$, diffusion is not able to dissipate enough of the energy stemming from the smooth forcing in this unidimensional setup. Actually, using only the Green function of the Laplacian, it can be shown that the variance will increase as fast as $\sqrt{t}$ as time $t$ goes on.
\end{remark}

\section{Dynamical Fractional Fields}
 \label{Sec:LinearFractDyn}

We have presented in Section \ref{Sec:HamiltonianL} a mechanism, governed by a linear operator $\mathcal L$ (Eq. \ref{eq:SecDefL}), able to transfer energy from the large scale $L$, a characteristic scale of the forcing $f$, towards the small ones. In the inviscid case ($\nu=0$), this dynamics is Hamiltonian when the force is shut down, and once forced generates a solution $u_0(t,x)$ that shares several properties with the white noise, as they are listed in Proposition \ref{eq:PDEHamilL}. This section is devoted to the presentation of the action of the linear operator $P_H$ (Eq. \ref{eq:IntroDefPH}) on this particular solution $u_0(t,x)$, recalling here its expression for any function $u$,
\begin{equation}\label{eq:SecDefPH}
P_Hu(t,x) \equiv \int e^{2i\pi k x}\frac{1}{|k|_{1/L}^{H+1/2}}\hat{u}(t,k)dk,
\end{equation}
where is introduced a regularized norm $|\cdot|_{1/L}$ of $k$ over the characteristic length of the forcing $L$. We do not need to precise the exact expression of this regularized norm in subsequent calculations, but require it to behave as the proper norm $|k|$ at large arguments, and that it goes to a finite positive value of order $1/L$ as $|k|$ goes to zero. To set ideas, we can keep in mind the expression $|k|_{1/L}^2=|k|^2+1/L^2$, a regularization that we will make use of in forthcoming numerical simulations. As we will see, the action of the operator $P_H$ on $u_0(t,x)$ will eventually generate a fractional Gaussian field at large times, and the power-law decrease that enters its Fourier transform will govern the regularity of this field at small scales. To show this, we consider in the following paragraph such a field, and derive for the sake of completeness its statistical properties.

\subsection{Fractional Gaussian fields}\label{Par:fGf}

\begin{proposition}\label{Prop:fGf}\textit{(Concerning an imaginary fractional Gaussian field of parameter $H$)} Consider the following Gaussian field
\begin{align}\label{eq:DefVFGF}
v_{H}(t,x) = (P_Hu_0)(t,x)\equiv \int e^{2i\pi kx}\frac{1}{|k|_{1/L}^{H+1/2}}\hat{u}_0(t,k)dk,
\end{align}
where $u_0(t,x)$ is the unique solution of the SPDE (Eq. \ref{eq:PDEHamilL}) starting with vanishing initial condition and sustained by a forcing term $f$. The field $v_H$ being defined as a linear operation on a Gaussian field $u_0$, is itself Gaussian, statistically homogeneous and of zero average. It is also a finite variance process for any $H>0$ and at any time, its value is given asymptotically by
\begin{align}\label{eq:CalcVarFGF}
\lim_{t\to\infty}\E|v_H(t,x)|^2 =  \frac{\mathcal C_f(0)}{2|c|}\int \frac{1}{|k|_{1/L}^{2H+1}}dk<+\infty.
\end{align}
Furthermore, the field $v_H(t,x)$ has locally in space the same regularity as a fractional Brownian motion of parameter $H$. To see this, define the increment over the scale $\ell$ as $\delta_\ell v_H(t,x)=v_H(t,x+\ell)-v_H(t,x)$, and get
\begin{align}
\delta_\ell v_H(t,x) =  \int e^{2i\pi kx}\frac{e^{2i\pi k\ell}-1}{|k|_{1/L}^{H+1/2}}\hat{u}_0(t,k)dk.
\end{align}
We have, for $H\in]0,1[$, the following behavior at small scales:
\begin{align}\label{eq:CalcVarIncrFGF}
\lim_{t\to\infty}\E|\delta_\ell v_H(t,x)|^2 \build{\sim}_{\ell\to 0}^{} |\ell|^{2H}\frac{\mathcal C_f(0)}{2|c|} \int \frac{\left|e^{2i\pi k}-1\right|^2}{|k|^{2H+1}}dk,
\end{align}
independently of the precise form of the regularization at vanishing wavelengths.
\end{proposition}

\begin{remark} The proofs are again straightforward since the Gaussian field $v_H$ (Eq. \ref{eq:DefVFGF}) is defined as a linear operation on $u_0(t,x)$ which behaves as a white noise at large time. The calculation of the variance (Eq. \ref{eq:CalcVarFGF}) is the consequence of the behavior of its correlation function $\mathcal C_{u_0}(t,x)$ at large time (Eq. \ref{eq:BehavCuWN}). This expression makes sense since integrability is warranted at infinity by $H>0$ and at the origin because the power-law is regularized over $1/L$. Let us underline that the variance depends explicitly on the precise choice of the regularization procedure.
\end{remark}

\begin{remark}\label{rem:UniversfGf} Very similarly to the calculation of the variance, the increment variance is a consequence of Eq. \ref{eq:BehavCuWN}. The equivalent at small scales (Eq. \ref{eq:CalcVarIncrFGF}) can be easily obtained when rescaling the dummy variable $k$ by $1/\ell$. Notice that the integrability at the origin furthermore requires that $H<1$ and does not necessitate a regularization procedure. In this sense, we can say that the behavior at small scales is independent of the mechanism of regularization at large scales (i.e. small wavelengths $k$).
\end{remark}

\begin{remark} As we can see the fractional Gaussian field $v_H$ (Eq. \ref{eq:DefVFGF}) is bounded, of finite variance (Eq. \ref{eq:CalcVarFGF}) and nowhere differentiable. Instead it shares the same local regularity as a fractional Brownian motion of parameter $H$, as is pinpointed by the behavior at small scales of the second-order structure function (Eq. \ref{eq:CalcVarIncrFGF}). It reproduces in this sense the regularity of a turbulent velocity field (Eq. \ref{eq:S2Turbu}) if we choose the particular value $H=1/3$. Because it is Gaussian, higher order structure functions $\E|\delta_\ell v_H|^q$ behave similarly as $(\E|\delta_\ell v_H|^2)^{q/2}$, which implies a spectrum of exponents $\zeta_q$ (Eq. \ref{eq:SqTurbu}) that depends linearly on $q$, at odds with experimental observations of turbulence.
\end{remark}

\begin{remark} \label{Rem:fGfH0} The boundary case $H=0$ is worth being considered. It enters in the construction of multifractal fields, as we will see in the next Section. Whereas the variance of such a field remains finite as $t\to \infty$ for $H>0$ (Eq. \ref{eq:CalcVarFGF}), it is no more the case for $H=0$. Instead, we get the diverging behavior
\begin{align}\label{eq:CalcVarFGF0DivTime}
\E|v_0(t,x)|^2 \build{\sim}_{t\to\infty}^{}  \frac{\mathcal C_f(0)}{|c|}\ln(ct),
\end{align}
and thus asymptotically as time gets large, $v_0(t,x)$ should be seen as a random distribution. Nevertheless, the covariance makes sense as a function away from the origin and reads, for $\ell\ne 0$,
\begin{align}\label{eq:CovarVFGFH0}
\lim_{t\to\infty}\mathcal C_{v_0}(t,\ell) \equiv\lim_{t\to\infty}\E[v_0(t,0)v_0^*(t,\ell)]= \frac{\mathcal C_f(0)}{2|c|}\int e^{2i\pi k\ell}\frac{1}{|k|_{1/L}}dk.
\end{align}
To see the infinite value of the variance in this limit (Eq. \ref{eq:CalcVarFGF0DivTime}), rescale the dummy variable $k$ by $1/\ell$ in Eq. \ref{eq:CovarVFGFH0}, and remark that the integral is then governed by the behavior of $|k|^{-1}_{\ell/L}$ near the origin when $\ell\to 0$, such that
\begin{align}\label{eq:CovarVFGFH0Near0}
\lim_{t\to\infty}\mathcal C_{v_0}(t,\ell) \build{\sim}_{\ell\to 0}^{} \frac{\mathcal C_f(0)}{|c|}\ln\left(\frac{L}{|\ell|}\right).
\end{align}
Thus, for $H=0$ and in the limit $t\to\infty$, the fractional Gaussian field (Eq. \ref{eq:DefVFGF}) has an infinite variance and is logarithmically correlated. In the following, we will find it convenient to write this asymptotic limit as
\begin{align}\label{eq:CovarVFGFH0DefH}
\lim_{t\to\infty}\mathcal C_{v_0}(t,\ell) = \frac{\mathcal C_f(0)}{|c|}\ln_+\left(\frac{L}{|\ell|}\right) +h(\ell),
\end{align}
where a smoothly-truncated logarithmic function $\ln_+(|x|)$ is introduced, which behaves as $\ln(|x|)$ as $|x|\to\infty$ and smoothly goes to zero as $|x|\to 0$, and $h(x)$ is a bounded and even function of its argument at least twice differentiable. As we will see, the very shapes of the truncation of the logarithm and of the function $h$ are not important, although they could be derived from Eq. \ref{eq:CovarVFGFH0}. Only the values at the origin of $h$ and its derivatives will matter, and we can get their exact expressions using a symbolic calculation software.
\end{remark}

\begin{remark}
We could have alternatively considered the field $\tilde{v}_H(t,x)$ defined in a similar manner as $v_H(t,x)$ (Eq. \ref{eq:DefVFGF}) but with an odd version $\tilde{P}_H$ of the operator $P_H$, which reads
\begin{align}\label{eq:DefVFGFHermitic}
\tilde{v}_{H}(t,x) = (\tilde{P}_Hu_0)(t,x)\equiv -i\int e^{2i\pi kx}\frac{k}{|k|_{1/L}^{H+3/2}}\hat{u}_0(t,k)dk,
\end{align}
without changing the global picture provided in Proposition \ref{Prop:fGf}. In particular, the variance is finite for $H>0$ and can be expressed similarly as in Eq. \ref{eq:CalcVarFGF} while slightly modifying the multiplicative factor, with the behavior of the second-order structure function being unchanged (Eq. \ref{eq:CalcVarIncrFGF}) for $H\in]0,1[$. When $H=0$, which corresponds to considering the action of the respective operator $\tilde{P}_0$ that we have defined in Eq. \ref{eq:IntroDefP0tilde}, we obtain the same logarithmic behaviors observed on the variance (Eq. \ref{eq:CalcVarFGF0DivTime}) and the correlation function (Eq. \ref{eq:CovarVFGFH0Near0}). Only the very shape of the additional function $h$ entering in the asymptotic limiting behavior of the respective correlation function $\mathcal C_{v_0}(t,\ell)$ (Eq. \ref{eq:CovarVFGFH0DefH}) is impacted by the possible parity of $P_H$, and we obtain instead that
\begin{align}\label{eq:CovarVFGFH0DefHTilde}
\lim_{t\to\infty}\mathcal C_{\tilde{v}_0}(t,\ell) = \frac{\mathcal C_f(0)}{|c|}\ln_+\left(\frac{L}{|\ell|}\right) +\tilde{h}(\ell),
\end{align}
where similarly $\tilde{h}(x)$ is a bounded and even function of its argument, at least twice differentiable.

\end{remark}

\subsection{Induced fractional dynamics}\label{Par:IndFracDynfGf}

In the light of the results of Section \ref{Sec:HamiltonianL} devoted to the design of a linear PDE, governed by the operator $\mathcal L$ (Eq. \ref{eq:SecDefL}), whose solution $u_0(t,x)$ behaves at large time as a white noise once it is forced, and the construction of fractional Gaussian fields (Paragraph \ref{Par:fGf}) using the linear action of the operator $P_H$ (Eq. \ref{eq:SecDefPH}) on $u_0(t,x)$, it is then tempting to consider the following dynamics
\begin{align}\label{eq:FracDynLPHnuF}
\partial_tu_{H,\nu}=P_H\mathcal LP_H^{-1}u_{H,\nu}+\nu\partial_x^2u_{H,\nu}+f.
\end{align}
Contrary to the Hamiltonian dynamics generated by the operator $\mathcal L$ in the inviscid and unforced situation, the operator $P_H\mathcal LP_H^{-1}$ does not preserve $|u_{H,0}(t,x)|^2$ in time. Actually, we will see that once forced, this dynamics will converge at large time towards a finite variance process, thus without the additional action of viscosity.

\begin{proposition}(Concerning the inviscid and forced fractional dynamics)\label{Prop:StatPropFracDyn} Consider the evolution
\begin{align}\label{eq:FracDynLPHNonuF}
\partial_tu_{H,0}=P_H\mathcal LP_H^{-1}u_{H,0}+f,
\end{align}
where $f(t,x)$ is a Gaussian random force defined in Eq. \ref{eq:CovGaussForce}, with $\mathcal C_f$ a real and even function of its argument, the linear operator $\mathcal L$ defined in Eq. \ref{eq:SecDefL}, and the operator $P_H$ defined in Eq. \ref{eq:SecDefPH}. For later convenience, once expressed in Fourier space, the dynamics provided in Eq. \ref{eq:FracDynLPHNonuF} reads
\begin{equation}\label{eq:FracDynEDPFourier}
\partial_t \hat{u}_{H,0}+c\partial_k\hat{u}_{H,0}+c\left(H+\frac{1}{2}\right)\frac{k}{|k|_{1/L}^2}\hat{u}_{H,0}=\hat{f}.
\end{equation}

Starting from the initial condition $u_{H,0}(0,x)=0$, the solution of this evolution is a zero-average Gaussian field $u_{H,0}(t,x)$.
Its correlation in space (Eq. \ref{eq:ExprCorrUXMY}) is a function of the difference of the positions only and is conveniently expressed in Fourier space as
\begin{align}\label{eq:ExprCorrUXFracCase}
\hat{\mathcal C}_{u_{H,0}}(t,k)=|k|_{1/L}^{-(2H+1)}\int_0^t|k-cs|_{1/L}^{2H+1}\hat{\mathcal C}_f(k-cs)ds.
\end{align}
Furthermore, for any $H>0$, the solution $u$ goes towards a statistically stationary regime as $t\to\infty$, such that its variance is finite, i.e.
\begin{align}\label{eq:FiniteVarFracCase}
\lim_{t\to\infty}\E|u_{H,0}|^2<\infty,
\end{align}
and its correlation function is given by
\begin{align}\label{eq:LimInfTimeExprCorrUXFracCase}
\lim_{t\to\infty}\hat{\mathcal C}_{u_{H,0}}(t,k)=\left\{
    \begin{array}{ll}
    \frac{1}{c}|k|_{1/L}^{-(2H+1)}\int_{-\infty}^k |s|_{1/L}^{2H+1}\hat{\mathcal C}_f(s)ds & \mbox{if } c>0 \\
     \frac{1}{-c}|k|_{1/L}^{-(2H+1)}\int_k^{\infty} |s|_{1/L}^{2H+1}\hat{\mathcal C}_f(s)ds & \mbox{if }c<0.
    \end{array}
\right.
\end{align}
Similarly to the fractional Gaussian field (see Proposition \ref{Prop:fGf} and Remark \ref{rem:UniversfGf}), as $t\to\infty$ and for $H\in]0,1[$, the solution $u$ shares the same local regularity as a fractional Brownian motion of parameter $H$, independently of the precise regularization procedure taking place at the scale $L$. Consequently, whatever the sign of $c$, the second order structure function behaves at small scales as
\begin{align}\label{eq:S2FracCaseInviscid}
\lim_{t\to\infty}\E|\delta_\ell u_{H,0}|^2&\build{\sim}_{\ell\to 0^+}^{} c_{H} \ell^{2H},
\end{align}
where the factor $c_H$ can be derived explicitly while introducing the function $\Gamma(z)=\int_0^\infty x^{z-1}e^{-x}dx$, and reads
\begin{align}\label{eq:MultiCstS2FracCaseInviscid}
c_{H} =\frac{1}{2|c|}\frac{(2\pi)^{2H+1}}{\sin(\pi H)\Gamma(1+2H)}\int_{-\infty}^{\infty}|s|_{1/L}^{2H+1}\widehat{\mathcal C}_f(s)\;ds.
\end{align}

\end{proposition}

\begin{remark} The proofs are once again straightforward since the PDE (Eq. \ref{eq:FracDynLPHNonuF}) is linear. Using the expression of the dynamics in Fourier space (Eq. \ref{eq:FracDynEDPFourier}), the unique solution $\hat{u}(t,k)$, starting with $\hat{u}(0,k)$=0, can be written as
\begin{equation}\label{eq:ExplSolFracDynEDPInviscid}
\hat{u}_{H,0}(t,k)=|k|_{1/L}^{-(H+1/2)}\int_0^t |k-c(t-s)|_{1/L}^{H+1/2}\hat{f}(s,k-c(t-s))\;ds.
\end{equation}
In a different setting, considering vanishing forcing but random initial conditions, it can be shown that $u_{H,0}$ remains bounded at any time.
The expression of the Fourier transform of the correlation function $\hat{\mathcal C}_{u_{H,0}}(t,k)$ (Eq. \ref{eq:ExprCorrUXFracCase}) can be obtained from the exact solution (Eq. \ref{eq:ExplSolFracDynEDPInviscid}), and its limiting value in the statistically stationary regime (Eq. \ref{eq:LimInfTimeExprCorrUXFracCase}) can be similarly justified. Let us notice that this expression is not an even function of the wavelength $k$ because of the complex nature of the setup, in particular when $c>0$, we have the asymptotical behavior
\begin{align}\label{eq:EquivLargeKStatSolFracInviscid}
\lim_{t\to \infty}\hat{\mathcal C}_{u_{H,0}}(t,k) \build{\sim}_{k\to\infty}^{} \frac{1}{c}k^{-(2H+1)}\int_{-\infty}^{\infty}|s|_{1/L}^{2H+1}\hat{\mathcal C}_f(s)\;ds,
\end{align}
as is expected for a fractional Gaussian field, whereas the decay for $k\to-\infty$ is much faster and completely governed by the forcing correlation function $\hat{\mathcal C}_f$. Similar behaviors are obtained when $c<0$ but looking at equivalents for large negative wavelengths. Integrability of $\hat{\mathcal C}_{u_{H,0}}(t,k)$ (Eq. \ref{eq:ExprCorrUXFracCase}) over $k\in\R$ requires $H>0$ and warrants a finite variance (Eq. \ref{eq:FiniteVarFracCase}). To compute the power-law behavior of the second-order structure function at small scales (Eq. \ref{eq:S2FracCaseInviscid}), including the multiplicative factor (Eq. \ref{eq:MultiCstS2FracCaseInviscid}), notice that at any time,
\begin{align}
\E|\delta_\ell u_{H,0}|^2=2   \int_{\R^+}\left[ 1- \cos\left( 2\pi k \ell\right)\right]\left[\hat{\mathcal C}_u(t,k)+\hat{\mathcal C}_u(t,-k)\right] dk,
\end{align}
and rescale the dummy variable $k$ by $\ell$ to obtain
\begin{align}\label{eq:ProofPowerLawS2FracInviscid}
\lim_{t\to\infty}\E|\delta_\ell u_{H,0}|^2&\build{\sim}_{\ell\to 0^+}^{} \ell^{2H}\frac{2}{|c|}\int_{-\infty}^{\infty}|s|_{1/L}^{2H+1}\widehat{\mathcal C}_f(s)\;ds   \int_{\R^+}\left[ 1- \cos\left( 2\pi k\right)\right]k^{-(2H+1)}dk,
\end{align}
where the last integral entering as a multiplicative factor in Eq. \ref{eq:ProofPowerLawS2FracInviscid} is finite for $H\in]0,1[$ and can be expressed with the help of the Gamma function $\Gamma$, which entails Eq. \ref{eq:MultiCstS2FracCaseInviscid}.
\end{remark}

\begin{remark}
Let us underline that, whereas the dynamics of Proposition \ref{Prop:HamDynL} (that generates white noise at large times) conserves energy, the corresponding underlying fractional dynamics of Proposition \ref{Prop:StatPropFracDyn}, without forcing, does not. This can be clearly seen from the expression of the latter in Fourier space (Eq. \ref{eq:FracDynEDPFourier}).
\end{remark}

\begin{remark}
As we can see, we succeeded in defining a Gaussian random field $u_{H,0}(t,x)$ as a solution of a PDE forced by a smooth term $f$ (Eq. \ref{eq:FracDynLPHNonuF}), that shares at large times several statistical properties with a fractional Gaussian field of parameter $H$ defined in Proposition \ref{Prop:fGf}, including a finite variance (Eq. \ref{eq:FiniteVarFracCase}) without the help of viscosity, and a local regularity governed by $H$ (Eq. \ref{eq:S2FracCaseInviscid}).
\end{remark}

\begin{remark} \label{Rem:AddNuFracDyn} It is then easy to include the effects of viscous diffusion while generalizing the results of Proposition \ref{eq:FracDynLPHNonuF}, considering the evolution provided in Eq. \ref{eq:FracDynLPHnuF}. Doing so, we introduce a new characteristic wavelength $k_\nu$  that gets larger as viscosity gets smaller, and is possibly dependent on $H$,  such that asymptotic behaviors as those of Eq. \ref{eq:LimInfTimeExprCorrUXFracCase} are expected in the finite range $1/L\ll |k|\ll k_\nu$. Equivalently using the terminology of turbulence phenomenology, the power-law behavior of the second-order structure function (Eq. \ref{eq:S2FracCaseInviscid}) is expected in the so-called inertial range $1/k_\nu\ll\ell\ll L$. Also, and for the same reasons, the statistically stationary regime can be reached at a finite time $t_\nu$ which depends on $\nu$. This will turn out to be very convenient from a numerical point of view.
\end{remark}

\begin{remark}
Moreover, working with a finite viscosity $\nu>0$ allows to revisit some key questions regarding turbulence phenomenology \cite{Fri95}. In particular, we could wonder what is the behavior of the viscous contribution  to the kinetic energy budget, i.e. $\nu\E |\partial_x u_{H,\nu}|^2$, as viscosity gets smaller. To do so, consider the solution of Eq. \ref{eq:FracDynLPHnuF}, starting form a vanishing initial condition, which reads
\begin{equation}\label{eq:ExplSolFracDynEDPNuPos}
\hat{u}_{H,\nu}(t,k)=\int_0^t \frac{|k-c(t-s)|_{1/L}^{H+1/2}}{|k|_{1/L}^{H+1/2}}e^{\frac{4\pi^2 \nu}{3c}\left[\left(k-c(t-s)\right)^3-k^3\right]}\hat{f}(s,k-c(t-s))\;ds.
\end{equation}
We notice that we recover Eq. \ref{eq:ExplSolFracDynEDPInviscid} from Eq. \ref{eq:ExplSolFracDynEDPNuPos} by considering the limit $\nu\to 0$. In this situation, the variance of the gradient $\partial_x u_{H,\nu}$ remains finite for any $H\in]0,1[$ and $\nu>0$ as $t\to \infty$, and we obtain
\begin{equation}\label{eq:ExplSolVarGradNuPos}
\lim_{t\to\infty}\E |\partial_x u_{H,\nu}|^2=4\pi^2\int_{k\in\R,s\in\R^+} (k+cs)^2\frac{|k|_{1/L}^{2H+1}}{|k+cs|_{1/L}^{2H+1}}e^{\frac{8\pi^2 \nu}{3c}\left[k^3-(k+cs)^3\right]}\hat{\mathcal C}_f(k)\;dsdk.
\end{equation}
To see how the gradient variance (Eq. \ref{eq:ExplSolVarGradNuPos}) diverges as $\nu\to 0$, rescale the dummy variable $s$ by $c\nu^{1/3}$ such that, for $H\in]0,1[$,
\begin{equation}\label{eq:EquivGradNuPosTo0}
\lim_{t\to\infty}\E |\partial_x u_{H,\nu}|^2\build{\sim}_{\nu\to 0}^{}\frac{4\pi^2}{c}\nu^{\frac{2}{3}(H-1)}\int_{s\in\R^+} s^{1-2H}e^{-\frac{8\pi^2 s^3}{3c}}\;ds\int_{k\in\R}|k|_{1/L}^{2H+1}\hat{\mathcal C}_f(k)\;dk.
\end{equation}
The equivalence provided in Eq. \ref{eq:EquivGradNuPosTo0} says that 
\begin{equation}\label{eq:DblLimGradNuPosTo0}
0<\lim_{\nu\to 0}\lim_{t\to\infty}\nu^{\frac{2}{3}(1-H)} \E |\partial_x u_{H,\nu}|^2 <\infty.
\end{equation}
Interestingly, the behavior of the viscous dissipation (Eq. \ref{eq:DblLimGradNuPosTo0}), choosing $H=1/3$ to get closer to a turbulent context, differs from the one expected from the Navier-Stokes equations \cite{Fri95}. In the latter, the gradient variance diverges as $\nu^{-1}$, whereas in the former, it diverges as $\nu^{\frac{2}{3}(H-1)}$. Regarding the kinetic energy budget of Eq.~\ref{eq:FracDynLPHnuF}, this also says that viscous dissipation does not participate as much as the contribution related to the fractional dynamics, in fact it participates less and less as $\nu$ gets smaller and smaller, in order to balance the energy injection provided by the force.
\end{remark}

\begin{remark}
Contrary to the regularity in space which is given by the parameter $H$, as it can be shown from the exact expression provided in Eq.~\eqref{eq:ExplSolFracDynEDPInviscid}, the regularity in time is governed by the rapid decorrelation of the forcing, and is thus independent on $H$. The temporal correlation structure of the solution is further investigated in Ref. \cite{ApoChev21}, where it is shown that, for a white-in-time forcing, a local regularity equal to that of Brownian motion is observed, while for correlated in time forcing a smooth in time solution is generated.
\end{remark}

\begin{remark}
For reasons that will become clear later, it will be useful to revisit the results listed in Proposition \ref{Prop:StatPropFracDyn} for the boundary case $H=0$, as the statistical properties of fractional Gaussian fields (Proposition \ref{Prop:fGf}) were revisited for this very particular value of the parameter $H$ (See Remark \ref{Rem:fGfH0}). In this case, at a finite time $t$, we can obtain the Fourier transform of the correlation function of $u_{0,0}(t,x)$, solution of Eq. \ref{eq:FracDynLPHNonuF}, using the expression provided in Eq. \ref{eq:ExprCorrUXFracCase}, and obtain
\begin{align}\label{eq:ExprCorrUXFracCaseH0}
\hat{\mathcal C}_{u_{0,0}}(t,k) =\frac{1}{c}|k|_{1/L}^{-1}\int_{k-ct}^k|s|_{1/L}\hat{\mathcal C}_f(s)ds,
\end{align}
which converges towards a bounded function as $t\to\infty$, depending on the sign of $c$, as was written in Eq. \ref{eq:LimInfTimeExprCorrUXFracCase}. Nonetheless, contrary to the case $H\in]0,1[$ for which the variance converges at large times towards a finite value (Eq. \ref{eq:FiniteVarFracCase}), the finiteness of the variance is no more warranted when $H=0$. Instead, we get, for any $|c|>0$, the following logarithmically diverging behavior at large times
\begin{align}\label{eq:DivVarFracCaseH0}
\E|u_{0,0}(t,x)|^2\build{\sim}_{t\to\infty}^{}\frac{1}{|c|}\ln(|c|t)\int_{\R}|s|_{1/L}\hat{\mathcal C}_f(s)ds.
\end{align}
Whereas the variance diverges with time (Eq. \ref{eq:DivVarFracCaseH0}), the correlation function (Eq. \ref{eq:ExprCorrUXMY}) is bounded away from the origin, and we can write, for any $\ell>0$,
\begin{align}\label{eq:ExprCorrUXPSFracCaseH0}
\lim_{t\to\infty}\mathcal C_{u_{0,0}}(t,\ell)=\int_{k\in\R}e^{2i\pi k\ell}\lim_{t\to\infty}\hat{\mathcal C}_{u_{0,0}}(t,k)dk.
\end{align}
We recover then the infinite value for the variance (Eq. \ref{eq:DivVarFracCaseH0}) in this regime while obtaining a logarithmically diverging behavior of the correlation function (Eq. \ref{eq:ExprCorrUXPSFracCaseH0}) at small scales, that is
\begin{align}\label{eq:DivergenceLCorrUXPSFracCaseH0}
\lim_{t\to\infty}\mathcal C_{u_{0,0}}(t,\ell)\build{\sim}_{|\ell|\to 0}^{}\frac{1}{|c|}\ln(L/|\ell|)\int_{\R}|s|_{1/L}\hat{\mathcal C}_f(s)ds.
\end{align}
We can thus see that the inviscid dynamics $\partial_t u_{0,0}=P_0\mathcal L P_0^{-1}u_{0,0}+f$ ultimately generates as time goes on a Gaussian field of infinite variance (Eq. \ref{eq:DivVarFracCaseH0}) which is logarithmically correlated in space (Eq. \ref{eq:DivergenceLCorrUXPSFracCaseH0}). It thus shares similar statistical properties with the fractional Gaussian fields (Proposition \ref{Prop:fGf}) of vanishing parameter $H=0$, as detailed in Remark \ref{Rem:fGfH0}.
\end{remark}

\section{Multifractal random fields and the induced nonlinear dynamics}
 \label{Sec:NonLinearMultiDyn}
This Section is devoted to the design of an additional nonlinear interaction term in the fractional evolution of Proposition \ref{Prop:StatPropFracDyn}  able to reproduce the observed multifractal nature of fluid turbulence. As mentioned earlier, the Gaussian framework that has been developed implies necessarily a  spectrum of exponents $\zeta_q$ of high-order structure functions (Eq. \ref{eq:SqTurbu}) that behaves linearly with the order $q$.
Inspired by the structure of probabilistic objects known as Multiplicative Chaos (MC) and related Multifractal processes, we will end up with a quadratic interaction that once added to the forced fractional linear evolution (Eq. \ref{eq:FracDynLPHNonuF}) supports the development of non-Gaussian fluctuations. As it is mentioned in the Introduction, MC is one of the first such probabilistic constructions which is statistically homogeneous. It was introduced in a turbulent context \cite{Man72} and further developed from a mathematical point of view \cite{Kah85,RhoVar14}, to reproduce the probabilistic nature of the turbulent energy dissipation as a field with lognormal statistics and a long range correlation structure. Again, roughly speaking, it is defined as the exponential of a Gaussian field with logarithmic covariance.
Numerical investigations detailed in the next Section furthermore indicate the multifractal behavior of this nonlinear evolution.

 \subsection{Complex Gaussian Multiplicative Chaos}

 \begin{proposition}\label{Prop:ComplexMC}(About a complex version of the Gaussian Multiplicative Chaos) Consider the following complex random field
 \begin{align}\label{eq:DefCGMC}
 M_\gamma(t,x) = e^{\gamma v_0(t,x)},
 \end{align}
 where $v_0(t,x)=P_0u_0(t,x)$ is a fractional Gaussian field of parameter $H=0$ (Eq. \ref{eq:DefVFGF}) whose asymptotic logarithmic correlation structure is detailed in Remark \ref{Rem:fGfH0}, and $\gamma\in\R$ an additional parameter.
As time goes to infinity, $M_\gamma(t,x)$ behaves as a random distribution such that, for any $\gamma\in\R$ and at any time and position,
\begin{align}\label{eq:AverageCGMC}
\E M_\gamma(t,x) = 1.
\end{align}
To see its distributional nature as $t\to\infty$, consider a $\mathcal C^{\infty}$ compactly supported function $g(x)$, of unit integral, and its rescaled version $g_\ell(x)\equiv g(x/\ell)/\ell$. We get the following asymptotic behavior, at large time $t\to\infty$ and at small scales $|\ell|\to 0$,  for any $q\in\N^*$ and $\gamma^2<\frac{|c|}{q\mathcal C_f(0)}$,
\begin{align}\label{eq:AsymptSmallScaleCGMC}
\lim_{t\to\infty}\E\left[\left|\int g_\ell(x) M_\gamma(t,x)dx \right|^{2q}\right] \build{\sim}_{|\ell|\to 0}^{}c_{q,\gamma}e^{q^2\gamma^2h(0)}\left(\frac{L}{|\ell|}\right)^{\frac{q^2\gamma^2\mathcal C_f(0)}{|c|}},
  \end{align}
where enters the value at the origin of the function $h$ defined in Eq. \ref{eq:CovarVFGFH0DefH} and a remaining multiplicative constant $c_{q,\gamma}$ given by
\begin{align}\label{eq:AsymptSmallScaleCGMCExpConst}
&c_{q,\gamma}=\\
&\int\prod_{i=1}^q\frac{1}{|x_i-y_i|^{\frac{\gamma^2\mathcal C_f(0)}{|c|}}} \prod_{i< j=1}^q\frac{1}{|x_i-y_j|^{\frac{\gamma^2\mathcal C_f(0)}{|c|}}}\frac{1}{|x_j-y_i|^{\frac{\gamma^2\mathcal C_f(0)}{|c|}}}\prod_{i=1}^{q}g(x_i)g(y_i) dx_idy_i.\notag
  \end{align}

\end{proposition}

\begin{remark} The proof of the statistical properties of the Complex Gaussian Multiplicative Chaos (GMC) $ M_\gamma(t,x)$ (Eq. \ref{eq:DefCGMC}) relies on several ingredients. Notice first that, as a consequence of the independence of the real and imaginary parts of $u_0(t,x)$ (Eq.  \ref{eq:ExprCorrUReIm}), entering in the definition of the fractional Gaussian field $v_0(t,x)$, the real and imaginary parts of $v_0(t,x)$ are statistically independent, and moreover, at any time and positions,\begin{align}\label{eq:IndepReImV0}
\E \left[v_0(t,x)v_0(t,y) \right]=0,
\end{align}
whereas we have formerly noted
\begin{align}\label{eq:RecallCorrV0}
\mathcal C_{v_0}(t,x-y)\equiv \E \left[v_0(t,x)v_0^*(t,y) \right]\in\C.
\end{align}
In the following, we will also rely on a particular statistical property of complex Gaussian random variables. Consider thus a complex Gaussian random variables $u$ such that $\E(u)=0$. We have the useful result
\begin{align}\label{eq:ExptExpCGV}
\E \left( e^{u}\right) = e^{\frac{1}{2}\E(u^2)},
\end{align}
which leads, at second order, to
\begin{align}
\E\left[\left|\int g_\ell(x) M_\gamma(t,x)dx \right|^{2}\right] &=\int g_\ell(x)g_\ell(y) \E\left[ e^{\gamma \left( v_0(t,x)+v_0^*(t,y)\right)}\right]dxdy\notag \\
&=\int g_\ell(x)g_\ell(y)e^{\gamma^2\mathcal C_{v_0}(t,x-y)}dxdy\notag,
  \end{align}
where we notice that $u=v_0(t,x)+v_0^*(t,y)$ is a zero average complex Gaussian random variable and used the properties of Eqs. \ref{eq:IndepReImV0} and \ref{eq:ExptExpCGV}. Relying then on the asymptotic form of $\mathcal C_{v_0}(t,x-y)$ at large time (Eq. \ref{eq:CovarVFGFH0DefH}), rescaling the dummy variables $x$ and $y$ by $\ell$, we obtain
\begin{align}\label{eq:SecondOrderCGMC}
\lim_{t\to\infty}\E\left[\left|\int g_\ell(x) M_\gamma(t,x)dx \right|^{2}\right] &\build{\sim}_{|\ell|\to 0}^{} \left(\frac{L}{|\ell|}\right)^{\frac{\gamma^2\mathcal C_f(0)}{|c|}}e^{\gamma^2h(0)}\int \frac{g(x)g(y) }{|x-y|^{\frac{\gamma^2\mathcal C_f(0)}{|c|}}}dxdy,
\end{align}
which requires $\gamma^2<|c|/\mathcal C_f(0)$ in order for the remaining integral to be finite.

Similar techniques can be used to derive higher-order moments. In particular we have
\begin{align}
\E&\left[\left|\int g_\ell(x) M_\gamma(t,x)dx \right|^{2q}\right] =\int \E\left[ e^{\gamma\sum_{i=1}^q v_0(t,x_i)+v_0^*(t,y_i)}\right] \prod_{i=1}^{q}g_\ell(x_i)g_\ell(y_i) dx_idy_i\notag\\
 &=\int e^{\frac{\gamma^2}{2}\sum_{i, j=1}^q\mathcal C_{v_0}(t,x_i-y_j)+\mathcal C_{v_0}(t,x_j-y_i)} \prod_{i=1}^{q}g_\ell(x_i)g_\ell(y_i) dx_idy_i\notag,
  \end{align}
such that
\begin{align}\label{eq:HighOrderCGMC}
&\lim_{t\to\infty}\E\left[\left|\int g_\ell(x) M_\gamma(t,x)dx \right|^{2q}\right] \build{\sim}_{|\ell|\to 0}^{} \left(\frac{L}{|\ell|}\right)^{\frac{q^2\gamma^2\mathcal C_f(0)}{|c|}}e^{q^2\gamma^2h(0)}\times \\
&\int\prod_{i=1}^q\frac{1}{|x_i-y_i|^{\frac{\gamma^2\mathcal C_f(0)}{|c|}}} \prod_{i< j=1}^q\frac{1}{|x_i-y_j|^{\frac{\gamma^2\mathcal C_f(0)}{|c|}}}\frac{1}{|x_j-y_i|^{\frac{\gamma^2\mathcal C_f(0)}{|c|}}}\prod_{i=1}^{q}g(x_i)g(y_i) dx_idy_i.     \notag
\end{align}
The finiteness of the remaining multiple integral entering on the RHS of Eq. \ref{eq:HighOrderCGMC}, and the implied range of possible values for the free parameter $\gamma$, is difficult to determine. Performing an integration over one variable, say $x_1$, and then making a spherical change of coordinates over the remaining $2q-1$ variables, would give $q^2\gamma^2<(2q-1)|c|/\mathcal C_f(0)$ stemming from the integration over the radial component, which is optimistic since integration over the $2q-2$ angles is not discussed. For $q\ge 3$, using a bound for the integrand which is simpler to analyse, as has been done in Lemma A.8 of Ref. \cite{GuMou16}, would instead give the more pessimistic range $q^2\gamma^2<q|c|/\mathcal C_f(0)$. A specially devoted communication on this would be needed, and is beyond the scope of the present article. This entails Eq. \ref{eq:AsymptSmallScaleCGMC} and motivates the proposed range of accessible values for $\gamma$.
\end{remark}

\begin{remark}
Similar behaviors as those depicted in Proposition \ref{Prop:ComplexMC} are again expected for a complex GMC $ \tilde{M}_\gamma$ based on the fractional Gaussian field  $\tilde{v}_0$ defined in Eq. \ref{eq:DefVFGFHermitic} for the particular value $H=0$, and which would read
\begin{align}\label{eq:DefCGMCtilde}
  \tilde{M}_\gamma(t,x) = e^{\gamma  \tilde{v}_0(t,x)}.
 \end{align}
Its distributional nature as $t\to\infty$ would also be characterized by Eq. \ref{eq:AsymptSmallScaleCGMC}, with a power-law exponent having the same quadratic dependence on the order $q$, and the same multiplicative constant  $c_{q,\gamma}$ (Eq. \ref{eq:AsymptSmallScaleCGMCExpConst}). Only the very shape of the function $h$ entering in Eq. \ref{eq:AsymptSmallScaleCGMC} and defined in Eq. \ref{eq:CovarVFGFH0DefH} would be impacted, requiring the use of the function $\tilde{h}$ defined in Eq. \ref{eq:CovarVFGFH0DefHTilde}.
\end{remark}

 \subsection{Construction of a complex multifractal field, and the calculation of its statistical properties} \label{Par:ComplexMC}

 \begin{proposition}\label{Def:MultiFieldV}\textit{(About a complex multifractal process $v_{H,\gamma}$)} Consider the random field $v_{H,\gamma}(t,x)$, defined as
\begin{align}
v_{H,\gamma}(t,x)&=P_H\left(e^{\gamma \tilde{P}_0u_0}u_0\right)(t,x)\equiv  \int e^{2i\pi kx}\frac{1}{|k|_{1/L}^{H+1/2}}\mathcal F\left[e^{\gamma \tilde{P}_0u_0} u_0\right](t,k)dk\notag\\
&=\int P_H(x-y)e^{\gamma \tilde{P}_0u_0(t,y)} u_0(t,y)dy,\label{eq:VHGammaDef}
\end{align}
where we have introduced the abusive notation
$$ P_H(x)=\int e^{2i\pi kx}\frac{1}{|k|_{1/L}^{H+1/2}}dk.$$
The field $v_{H,\gamma}(t,x)$ is statistically homogeneous and its average is 0 at any time $t$. As time gets large, the field is such that, for any $q\in\N^*$, $H\in]0,1[$ and $\gamma^2\frac{\mathcal C_f(0)}{|c|}<\min(2H/q,1)$,
\begin{align}\label{eq:FinitMomentsProposition}
\lim_{t\to\infty}\E\left|v_{H,\gamma}(t,x)\right|^{2q}<\infty.
\end{align}
Furthermore, for the same range of values of the parameters $H$ and $\gamma$, the random field $v_{H,\gamma}$ (Eq. \ref{eq:VHGammaDef}) exhibits a multifractal local regularity, as can be quantified by its respective structure functions of order $2q$, which behave at small scales as
\begin{align}\label{eq:HOStatVHGamma}
\lim_{t\to\infty}\E\left|\delta_\ell v_{H,\gamma}\right|^{2q}\build{\sim}_{\ell\to 0}^{}c_{H,\gamma,q}\ell^{2qH}\left( \frac{\ell}{L}\right)^{-q^2\gamma^2\frac{\mathcal C_f(0)}{|c|}},
\end{align}
where the explicit expression of the multiplicative constant $c_{H,\gamma,q}$ at the second order $q=1$ is provided in Eq. \ref{eq:EquivSmallLLimVarDeltaLVHgamma}.

Concerning the skewed nature of the probability laws of the random field $v_{H,\gamma}$ (Eq. \ref{eq:VHGammaDef}), we obtain in a similar way, again for $H\in]0,1[$, but for $\gamma^2\frac{\mathcal C_f(0)}{|c|}<\min(3H/2,1)$,  non trivial odd-order statistics, as they can be quantified by the behavior at small scales of the following expectation
\begin{align}\label{eq:OddStatVHGamma}
\lim_{t\to\infty}\E\left[\delta_\ell v_{H,\gamma}\left|\delta_\ell v_{H,\gamma}\right|^{2}\right]\build{\sim}_{\ell\to 0}^{}d_{H,\gamma}\ell^{3H}\left( \frac{\ell}{L}\right)^{-2\gamma^2\frac{\mathcal C_f(0)}{|c|}},
\end{align}
where $d_{H,\gamma}\in\R$ is a real and finite multiplicative factor, whose exact expression is given in Eq. \ref{eq:CalcDHGamma3rd}.
\end{proposition}

\begin{remark}\label{eq:LimProofsVHGamma}
We partially prove the statements of Proposition \ref{Def:MultiFieldV} in Appendix \ref{App:ProofVHGamma} while deriving the expressions of the expectations using the Gaussian integration by parts formula (see Lemma \ref {Lemma:GIPComplex}, which is adapted from Lemma 2.1 of Ref. \cite{RobVar08} for the general case of complex Gaussian random variables). In particular, we derive in an exact fashion the variance and increments variance, i.e. Eqs. \ref{eq:FinitMomentsProposition} and \ref{eq:HOStatVHGamma} using the particular value $q=1$.
Also, we justify in Appendix \ref{App:ProofVHGamma} the scaling behavior of the third-order structure function (Eq. \ref{eq:OddStatVHGamma}), computing in particular the value of the multiplicative factor  $d_{H,\gamma}$ and showing that it is finite for the proposed range of parameters $H$ and $\gamma$. Nonetheless, its expression is intricate. Whereas we demonstrate that its value is finite, we fail at giving simple arguments to justify that it does not vanish, also its sign remains unknown. For $q\ge 2$, expressions of moments (Eqs. \ref{eq:FinitMomentsProposition}) and structure functions  (Eq. \ref{eq:HOStatVHGamma}) get even more cumbersome. For this reason, we only provide heuristics that led us to propose the scaling behavior of Eq. \ref{eq:HOStatVHGamma}, and the range of values of $\gamma$ for which this asymptotic behavior is expected to make sense.
\end{remark}

\begin{remark}
As stated in Proposition \ref{Def:MultiFieldV}, given the limitations listed in Remark \ref{eq:LimProofsVHGamma}, the complex random field $v_{H,\gamma}$ (Eq. \ref{eq:VHGammaDef}) behaves at infinite time as a multifractal function (Eq. \ref{eq:HOStatVHGamma}), and furthermore its real part is skewed (Eq. \ref{eq:OddStatVHGamma}). Interestingly, as argued in Ref. \cite{CheGar19}, $v_{H,\gamma}$ has no natural equivalent in a purely real setup. Instead, defining a real, skewed and multifractal unidimensional random field requires a more sophisticated method of construction, as is developed in  Ref. \cite{CheGar19}. In this case, the real part $\Re v_{H,\gamma}$ of $v_{H,\gamma}$ can be seen as a realistic probabilistic representation of the longitudinal component of the three-dimensional turbulent velocity vector field, whose statistical properties are listed in Refs. \cite{Fri95,CheGar19} if the particular value $H=1/3+\gamma^2\frac{2\mathcal C_f(0)}{3|c|}$ is chosen.  In this case, the third-order structure function (Eq. \ref{eq:OddStatVHGamma}) behaves linearly with the scale $\ell$, as is suggested by the so-called four-fifths law of turbulence (see Ref. \cite{Fri95}), which governs the energy transfers through scales. In a turbulent setting, focusing again on the longitudinal component of the velocity field, as is measured in wind tunnels, the intermittency parameter is observed to be universal, i.e. independent of the Reynolds, and small, of the order of $2\gamma^2\frac{\mathcal C_f(0)}{|c|}\approx 0.025$ (see for instance Ref. \cite{CheCas12}). \end{remark}

\begin{remark}
Let us finally remark that even and odd-order statistics behave in a different manner, as was already observed in different, although similar, random fields \cite{CheGar19}. In particular, the third-order structure function goes towards 0 as $\ell^{3H-2\gamma^2\frac{\mathcal C_f(0)}{|c|}}$ (Eq. \ref{eq:OddStatVHGamma}), whereas it is expected heuristically that $\lim_{t\to\infty}\E\left|\delta_\ell v_{H,\gamma}\right|^{3}$ would go to zero as $\ell^{3H-\frac{9}{4}\gamma^2\frac{\mathcal C_f(0)}{|c|}}$. To see this, assume that Eq. \ref{eq:HOStatVHGamma} can be extended to non-integer values, and take $q=3/2$. Although surprising, this remark is consistent with the constraint that, at any time and any scale, we should have $\E\left[\Re \delta_\ell v_{H,\gamma}\left|\delta_\ell v_{H,\gamma}\right|^{2}\right] \le \E\left[\left|\delta_\ell v_{H,\gamma}\right|^{3}\right]$. Nonetheless, it remains to be rigorously shown, following a more general approach as the one developed in Ref. \cite{CheGar19}, that would give access to the behavior at small scales of $\lim_{t\to\infty}\E\left|\delta_\ell v_{H,\gamma}\right|^{q}$ for any $q\in\R$.
We keep this important  perspective for future investigations. We also draw attention to the fact that the right-hand side of Eq.~\eqref{eq:OddStatVHGamma} is real, while $v_{H,\gamma}$ is complex, eliciting the fact that only the real part of the increment $\delta_{\ell} v_{H,\gamma}$ is skewed, while its imaginary part is not.
\end{remark}

\begin{remark}
We notice that the structure function exponents of the multifractal field $v_{H,\gamma}$ do not depend on the precise shape of the correlation function of the external forcing, but only on its variance, given by $\mathcal C_f(0)$, which can be incorporated into the free parameter $\gamma$.  We recall that, in the phenomenology of turbulence, such exponents are expected to be independent on the properties of the forcing. This would thus correspond to choosing this free parameter $\gamma$ of the present model in units of the forcing variance.
\end{remark}

 \subsection{Induced nonlinear dynamics} \label{Par:InducedDynComplexMC}

Let us now explore the consequence of the probabilistic ansatz $v_{H,\gamma}(t,x)$ (Eq. \ref{eq:VHGammaDef}) that we recall to exhibit multifractal statistics (Eq. \ref{eq:HOStatVHGamma}), as presented in Proposition \ref{Def:MultiFieldV}. In particular, in the same way as we built up the fractional inviscid evolution in Eq. \ref{eq:FracDynLPHNonuF}, we would like to extract, at least heuristically, a dynamics for a field $u_{H,\gamma,0}(t,x)$ which once forced by $f$ would lead to similar statistics as $v_{H,\gamma}(t,x)$. Recall first that the field $u_0$ entering in the definition of $v_{H,\gamma}(t,x)$ (Eq. \ref{eq:VHGammaDef}) evolves in the inviscid and unforced situation according to
$$ \partial_t u_0=\mathcal Lu_0,$$
where the transfer operator $\mathcal L$ is defined in Eq. \ref{eq:SecDefL}. Accordingly, we thus expect that
\begin{align}
\partial_t \left( e^{\gamma \tilde{P}_0u_0} u_0\right) &= e^{\gamma \tilde{P}_0u_0} \partial_t u_0+\gamma \tilde{P}_0\partial_tu_0\left( e^{\gamma \tilde{P}_0u_0} u_0\right)\notag\\
&= \mathcal L\left( e^{\gamma \tilde{P}_0u_0} u_0\right)+\gamma \tilde{P}_0\mathcal L u_0\left( e^{\gamma \tilde{P}_0u_0} u_0\right).\notag
\end{align}

From a formal point of view, note $\mathcal W$ a functional of some complex function $h:\R\rightarrow\C$, implicitly defined as  $\mathcal W[h](x)e^{\gamma \tilde{P}_0\mathcal W[h](x)}=h(x)$, if it exists, such that,
\begin{align}\label{eq:NLPDEwithW}
\partial_tv_{H,\gamma}\equiv P_H \partial_t \left( e^{\gamma \tilde{P}_0u_0} u_0\right) &= P_H\mathcal LP_H^{-1}v_{H,\gamma}\\
&+\gamma P_H \left[\left( \tilde{P}_0\mathcal L \mathcal W\left[P_H^{-1}v_{H,\gamma}\right]\right)\left( P_H^{-1}v_{H,\gamma}\right)\right].\notag
\end{align}
We can see that the first term at the RHS of Eq. \ref{eq:NLPDEwithW}, which is linear in the variable $v_{H,\gamma}$  coincides with the deterministic part of the fractional evolution proposed in Eq. \ref{eq:FracDynLPHNonuF}. The second term, proportional to the multifractal parameter $\gamma$, is not clearly closed in terms of $v_{H,\gamma}$, but certainly introduces a nonlinearity in the picture.

 \subsection{A closure approach} \label{Par:ClosureInducedDynComplexMC}
The functional $\mathcal  W$ entering in the nonlinear evolution of Eq. \ref{eq:NLPDEwithW} resembles a functional generalization of the Lambert W-function, which is a multivalued function of $\C\rightarrow\C$. Much care is needed to make sense of it, and this is out of the scope of the present article. Although we could write formally the functional $\mathcal  W[h]$ in a recursive manner, a tractable form as a function of $h(x)$ is not known. Consequently, the evolution given in Eq. \ref{eq:NLPDEwithW} is not closed in terms of the field $v_{H,\gamma}$, and in order to close it, we propose to use the simplest and natural approximation given by
\begin{align}\label{eq:ApproxW}
\mathcal  W[h](x)\approx h(x),
\end{align}
which corresponds to making the approximation
\begin{align}\label{eq:SeriesExpW}
P_H^{-1}v_{H,\gamma}(t,x)=e^{\gamma \tilde{P}_0u_0(t,x)}u_0(t,x)\approx u_0(t,x).
\end{align}
In other words, we approximate the functional $\mathcal  W$ entering in Eq. \ref{eq:NLPDEwithW} by the identity (Eq. \ref{eq:ApproxW}), and this can be motivated by the Taylor series of the exponential entering in Eq. \ref{eq:SeriesExpW}, keeping only the first term in its development as powers of $\gamma$. Doing so, we end up with the following, approximate but closed, nonlinear evolution for the multifractal field $v_{H,\gamma}$ (Eq. \ref{eq:NLPDEwithW})
\begin{align}\label{eq:ApproxNLPDEVHGamma}
\partial_tv_{H,\gamma} \approx P_H\mathcal LP_H^{-1}v_{H,\gamma}+\gamma P_H \left[\left( \tilde{P}_0\mathcal LP_H^{-1}v_{H,\gamma}\right)\left( P_H^{-1}v_{H,\gamma}\right)\right].
\end{align}
The approximative evolution of the multifractal field, as given in Eq. \ref{eq:ApproxNLPDEVHGamma}, motivated us to propose the nonlinear PDE of the introduction, Eq. \ref{eq:IntroPDE}, forced by $f$ with the additional action of viscosity $\nu$. Of course, in presence of such an additional quadratic interaction, the theoretical analysis becomes much more demanding than the linear fractional evolution of Section \ref{Par:IndFracDynfGf}. This is why we will focus in the sequel on a numerical exploration of this stochastically forced nonlinear PDE.

\section{Numerical simulations}
 \label{Sec:Nums}

\subsection{Numerical setup}

The aim of this section is to present a numerical investigation of the statistical properties of the solution
$u_{H,\gamma,\nu}$ of Eq. \ref{eq:IntroPDE}. To do so, we make the following numerical proposition.

\begin{Numproposition} \label{Numprop:PDE}For periodic boundary conditions, over the period $L_{{\text{tot}}}$ of the spatial numerical domain, starting from the initial condition $u_{H,\gamma,\nu}(0,x)=0$, over the grid
$x=\{-N/2+1,\dots,0,\dots,N/2\}\Delta x$, where $N$ is the number of collocation points and $\Delta x=L_{{\text{tot}}}/N$, we are solving the numerical problem
\begin{align}\label{eq:PropNumPDE}
du_{H,\gamma,\nu}=&\left[P_H\mathcal LP_H^{-1}u_{H,\gamma,\nu}+ \gamma P_H\left[\left(\tilde{P}_0\mathcal L P_H^{-1}u_{H,\gamma,\nu}\right) \left( P_H^{-1}u_{H,\gamma,\nu}\right) \right] +\nu \partial^2_x u_{H,\gamma,\nu}\right]\Delta t \\
&+ f_{{\text{trunc}}}\sqrt{\Delta t},\notag
\end{align}
where the operators $\mathcal L$, $P_H$ and its inverse $P_H^{-1}$, and $\tilde{P}_0$ are defined respectively in Eqs. \ref{eq:IntroDefL}, \ref{eq:IntroDefPH}, \ref{eq:IntroDefPHinverse} and \ref{eq:IntroDefP0tilde}. The force $ f_{{\text{trunc}}}$ that sustains the dynamics is a truncated version of the forcing term $f(t,x)$ defined in Eq. \ref{eq:NumDefCf} which vanishes at the boundaries, i.e. $ f_{{\text{trunc}}}(t,\pm L_{{\text{tot}}}/2)=0$. The time marching is based on an explicit predictor-corrector algorithm in which a single instance of the force $f_{{\text{trunc}}}$ is generated at every time step $\Delta t$.
\end{Numproposition}

\begin{remark}
The evolution given in Eq. \ref{eq:PropNumPDE} involves several operations that are nonlocal in physical space, but local in Fourier space, including the convolutions with the operators $P_H$, its inverse, $\tilde{P}_0$ and the second derivative associated to viscous diffusion (recall that the Fourier symbol of the second derivative in physical space is $\mathcal F[\partial_x^2](k)=-(2\pi)^2k^2$). In this periodic framework, we will massively rely on the Discrete Fourier Transform (DFT) to evaluate the deterministic part at the RHS of Eq. \ref{eq:PropNumPDE}. Corresponding available wavelengths are thus given by $k=\{-N/2+1,\dots,0,\dots,N/2\}\Delta k$ where $\Delta k=1/L_{{\text{tot}}}$. For full benefit of the Fast Fourier Transform (FFT) algorithm to evaluate the DFT, we choose $N$ to be a power of 2, i.e. $N=2^n$. Evaluations of the transfer operator $\mathcal L$ and the quadratic term proportional to the parameter $\gamma$ are performed in physical space, which implies several back and forth computations of the DFT and its inverse, in what is known as a pseudo-spectral method. Furthermore, to get rid of the aliasing error induced by the quadratic nonlinear term, we use a de-aliasing procedure based on the $3/2$-rule (see for instance Ref. \cite{Pop00}).
\end{remark}

\begin{remark}
We also recall the definition of the forcing term $f(t,x)$ that enters in the continuous evolution of Eq. \ref{eq:IntroPDE}, which is a complex Gaussian random force defined in Eq. \ref{eq:CovGaussForce}. It is uncorrelated in time, each instance of the force is taken as
\begin{align}\label{eq:DefNumConvF}
f(t,x) = \int e^{-\frac{(x-y)^2}{2L^2}}dW(t,y),
\end{align}
where $dW(t,x)=\frac{1}{\sqrt{2}}\left[dW_r(t,x)+idW_i(t,x)\right]$ with $dW_r(t,x)$ and $dW_i(t,x)$ being at each time $t$ independent copies of the increment over $dx$ of a real Wiener process. Remark that another choice for the convolution kernel entering in Eq. \ref{eq:DefNumConvF} could be made, as long as its Fourier transform decreases rapidly above the characteristic wavelength $1/L$. Choosing the force $f$ as given in Eq. \ref{eq:DefNumConvF} implies that its real and imaginary parts are independent, or equivalently that $\E[f(t,x)f(t,y)]=0$ at any positions $x$ and $y$. From a numerical point of view, and in our periodic setup, an instance at a time $t$ of $f(t,x)$ is conveniently obtained by multiplying the DFTs of the convolution kernel and of $N$ independent instances of a zero-average Gaussian random variable $\mathcal N(0,\Delta x/2)$ of variance $\Delta x/2$ for both the real and imaginary parts. The form of the spatial dependence $\mathcal C_f$ of its covariance (Eq. \ref{eq:CovGaussForce}) is explicitly given by
\begin{align}\label{eq:NumDefCf}
\mathcal C_f(x)\equiv \sqrt{\pi L^2} e^{-\frac{x^2}{4L^2}},
\end{align}
and its expression in Fourier space corresponds to
\begin{align}\label{eq:NumDefHatCf}
\hat{\mathcal C}_f(k)=2\pi L^2e^{-4\pi^2k ^2L^2}.
\end{align}
\end{remark}

\begin{remark}
Notice that we have chosen the same characteristic large length scale $L$ in the definitions of the force $f$ (Eq. \ref{eq:NumDefCf}) and of the operator $P_H$ (Eq. \ref{eq:IntroDefPH}). It would be interesting to explore precisely the influence of choosing different scales to define forcing and fractional operators, although we expect from a physical point of view that these scales are of the same order.  With the particular choice of $L$ of a few fractions of $L_{{\text{tot}}}$, we are able to observe the cascading of energy towards small scales, as will be developed in the sequel. We believe that choosing $L$ or $L_{{\text{tot}}}$ in the definition of  $P_H$ (Eq. \ref{eq:IntroDefPH}) would give similar numerical results, as it can be fully shown in the linear framework (choose $\gamma=0$ in Eq. \ref{eq:PropNumPDE}) since statistical properties are in this case known in an exact fashion.
\end{remark}

\begin{remark}
Importantly, the periodization of the operator $\mathcal L$ (Eq. \ref{eq:IntroDefL}) introduces a discontinuity at the boundaries of the integration domain $x=\pm L_{{\text{tot}}}/2$. A first way to get rid of this spurious discontinuity is to consider a periodic version of this operator such as $\mathcal L_{{\text{per}}}=iL_{{\text{tot}}}\sin(2\pi x/L_{{\text{tot}}})$. Doing so, as was explored numerically in Ref. \cite{DyaZwo19} while solving an equation similar to the evolution given in Eq. \ref{eq:PDEHamilL}, the induced solution looses the convenient property of statistical homogeneity, and in particular energy accumulates at the boundaries. To prevent this accumulation of energy at the boundaries, while keeping an approximate statistically homogeneous region of space near the origin, i.e. far from the boundaries, and using the operator $\mathcal L$ as it is defined in Eq. \ref{eq:IntroDefL}, we propose to use, instead of the force $f(t,x)$ defined in Eq. \ref{eq:NumDefCf}, its truncated version
\begin{align}\label{eq:DefFtrunc}
f_{{\text{trunc}}}(t,x)=e^{-\frac{ x^2}{ L^2_{{\text{tot}}}/4-x^2}}f(t,x).
\end{align}
The particular choice of the bump function to truncate the force $f$ is not crucial at this stage. This choice is mostly motivated by the fact that $f_{{\text{trunc}}}(t,x)$ coincides with $f(t,x)$ at the origin, and goes smoothly towards zero at $x\to \pm L_{{\text{tot}}}$, without thus introducing another discontinuity at the boundaries. Using $f_{{\text{trunc}}}$ (Eq. \ref{eq:DefFtrunc}) instead of $f$ (Eq. \ref{eq:NumDefCf}) introduces nonetheless an inhomogeneous term in the evolution given in Eq. \ref{eq:PropNumPDE}. 

This inhomogeneity is nevertheless a feature of the numerical simulations only, the continuous solution over $x\in\R$ being statistically homogeneous, as it can be shown rigorously in the linear and Gaussian setup (see Proposition \ref{Prop:StatPropFracDyn}). Indeed, such a truncation is not required when dealing only with linear interactions, although the numerical solution eventually becomes  discontinuous at the boundaries of the periodical domain. Such a numerical framework has been alternatively chosen in Ref. \cite{ApoChev21}, in which the solution can be observed to be statistically homogeneous. Introducing such an inhomogeneous term $f_{{\text{trunc}}}$ (Eq. \ref{eq:DefFtrunc})  in the numerical problem furthermore guarantees numerical stability in the presence of an additional nonlinear term, as it is encountered when solving Eq.~\ref{eq:ApproxNLPDEVHGamma}. 

We will extensively comment in the next paragraph on the implications of taking a vanishing force $f_{{\text{trunc}}}$ at the boundaries of the numerical domain on the solution of the PDE under study (Eq. \ref{eq:PropNumPDE}). In particular, we will see that the solution will be observed in a good approximation to be homogeneous around the origin, say in the restricted domain $x\in[-0.2 L_{{\text{tot}}},0.2 L_{{\text{tot}}}]$. We will also observe that the solution of the PDE of Eq. \ref{eq:PropNumPDE}, with a vanishing initial condition, when sustained by the truncated force $f_{{\text{trunc}}}$ (Eq. \ref{eq:DefFtrunc}), will also vanish at the boundaries of the numerical domain $x=\pm L_{{\text{tot}}}/2$. This is guaranteed by the facts that the deterministic terms entering in the RHS of Eq. \ref{eq:PropNumPDE} must also vanish at the boundaries by periodicity and that the operator $\mathcal L$ is skew-symmetric.
\end{remark}

\begin{remark}
The addition of viscosity enforces that a stationary state is reached in a finite time, of the order of $1/c k_{\nu}$, where $k_{\nu}$ is the characteristic dissipative wavelength, both in the linear and nonlinear equations.
Even though we expect the forced Eq.~\ref{eq:ApproxNLPDEVHGamma} to reach a state of finite variance even in the absence of viscous dissipation (based on the multifractal ansatz of Eq.~\ref{eq:VHGammaDef}), this state can only be reached at infinite time, in order to populate in an appropriate manner infinitely small scales, as it happens in the linear case (Eq.~\ref{eq:FracDynLPHNonuF}). Furthermore, with viscosity not only the variance is finite but also, for instance, the second-order structure function at small scales.
\end{remark}

\subsection{Numerical results}

We perform several simulations of the numerical problem detailed in the proposition \ref{Numprop:PDE}. The parameters of the simulations are chosen in the following way. Without loss of generality, we take $L_{{\text{tot}}}=1$. The integral length scale is chosen as $L=L_{{\text{tot}}}/10$ and we consider the particular value $H=1/3$ as suggested by the phenomenology of turbulence. The rate of transfer of energy is set to $c=10$. Three different values for the intermittency parameter $\gamma=0$, $\sqrt{0.01}$ and $\sqrt{0.02}$ are chosen. Based on the numerical stability of the underlying heat equation when facing a discontinuity, the time step is expected to be chosen of the order of $\Delta t\approx(\Delta x)^2$, although this would require a prohibitive numerical cost. Instead, since the solution is expectedly continuous, we will choose $\Delta t=\Delta x$, a value that we found small enough to avoid singular behaviors and to be numerically tractable. Viscosity $\nu$ and number of collocation points $N$ are chosen such that the smallest length scale of the problem, which is of the order of $(\nu/|c|)^{1/3}$ (see the discussion provided in Remark \ref{Rem:AddNuFracDyn}) is properly resolved such that no numerical instabilities are observed. With the given choices made for the aforementioned parameters, we moreover consider the pairs of values $(N;\nu) = (2^{12};10^{-5})$, $(2^{13}; 10^{-6})$, $(2^{14};10^{-7})$, $(2^{16}; 10^{-8})$, $(2^{16}; 10^{-9})$. Notice that we have used the same resolution $N=2^{16}$ to study the values of the viscosity $\nu=10^{-8}$ and $10^{-9}$, because we observed for the former case some numerical instabilities. We checked that integration in time is long enough to reach a statistically steady regime, and only then various quantities of interest are averaged at several times such that the statistical samples are independent.

  \begin{figure}[t]
    \centering
 \includegraphics[width=\linewidth]{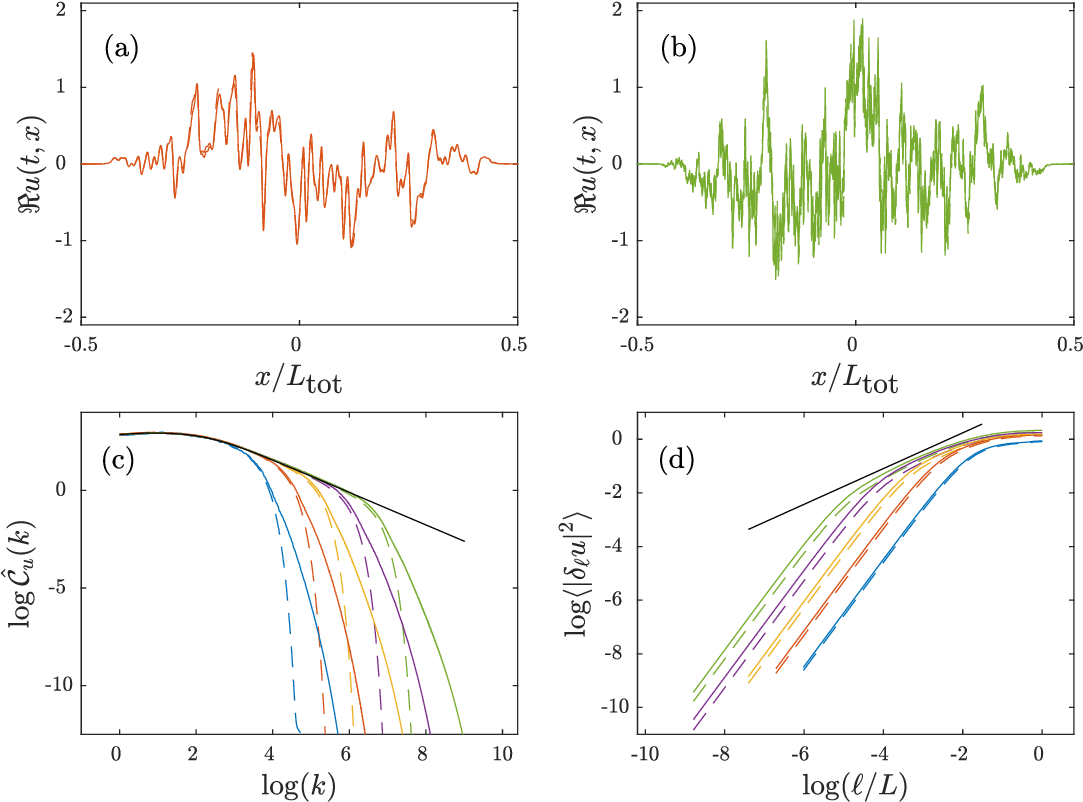}
    \caption{Local and statistical behaviors of the solution $ u_{H,\gamma,\nu}(t,x)$ detailed in Proposition \ref{Numprop:PDE}. (a): spatial profiles of  $\Re u_{H,\gamma,\nu}(t,x)$ at a given time $t$ in the statistically stationary regime for $L=L_{{\text{tot}}}/10$, $c=10$, $H=1/3$, $\nu=10^{-6}$, with $\gamma=0$ (dashed line)  and $\gamma=\sqrt{0.02}$ (solid line).  (b): same plot as in (a), but for a lower value of viscosity $\nu=10^{-9}$. (c): estimations of the power spectrum based on the averaged periodograms (see text) of the solution for various values of viscosity $\nu=10^{-5}$, $10^{-6}$, $10^{-7})$, $10^{-8}$ and $10^{-9}$ (from left to right), using dashed lines for $\gamma=0$, and a solid ones for  $\gamma=\sqrt{0.02}$. We superimpose with a solid black line the asymptotic prediction made in Eq. \ref{eq:LimInfTimeExprCorrUXFracCase}, which has been obtained in the Gaussian and fractional case, properly weighted by a multiplicative factor to take the truncation of the force into account (see text). (d) Similar plot as for (c) but for the second order structure function, i.e. the variance of the increments, following an averaging procedure detailed in the text. We superimpose the expected asymptotic power-law behavior, given in Eq. \ref{eq:S2FracCaseInviscid}, which is properly weighted (see text) and represented by a solid black line. }\label{fig:SecondOrder}
  \end{figure}

We display in Fig. \ref{fig:SecondOrder} the results of our simulations. We begin with the spatial representation of the solution in the statistically stationary regime at a given time $t$, at a moderate viscosity $\nu=10^{-6}$ (Fig. \ref{fig:SecondOrder}(a) in red) and for the lowest value $\nu=10^{-9}$ (Fig. \ref{fig:SecondOrder}(b) in green). For both cases, we moreover superimpose the Gaussian case $\gamma=0$ using a dashed-line and $\gamma=\sqrt{0.02}$ using a solid line. As we already explained, the solution $u_{H,\gamma,\nu}(t,x)$ vanishes at the boundaries $x=\pm L_{{\text{tot}}}/2$ of the domain. Also, we can barely see in this representation a difference between the Gaussian $\gamma=0$ and intermittent $\gamma\ne 0$ cases. Correspondingly, the dashed and solid lines almost perfectly superimpose. This shows from a numerical point of view that somehow the intermittent solution could be approached in a perturbative way with respect to the Gaussian solution. As viscosity decreases, we can also observe the appearance of fluctuations at smaller and smaller length scales, making overall the series of Fig. \ref{fig:SecondOrder}(b) rougher than those displayed in Fig. \ref{fig:SecondOrder}(a). Similar plots could be obtained for the imaginary parts $\Im u_{H,\gamma,\nu}(t,x)$ of the solution instead of the real one $\Re u_{H,\gamma,\nu}(t,x)$.

We present in Fig. \ref{fig:SecondOrder}(c) in a logarithmic representation the estimation of the power spectrum, i.e. the Fourier transform $\hat{\mathcal C}_{u_{H,\gamma,\nu}}(t,k)$ for various values of viscosity, and for both values $\gamma=0$ (dashed-line) and $\gamma=\sqrt{0.02}$ (solid-line). The estimation is made using the periodogram over the full periodical spatial domain, i.e. the square norm of the DFT of the solution, normalized by $L_{{\text{tot}}}$, which is averaged in time using independent instances. As viscosity decreases, a wider and wider range of energy-populated wavelengths develops, in a similar way as small-scale fluctuations appear in spatial profiles (Figs. \ref{fig:SecondOrder}(a) and (b)). For the smallest viscosity $\nu=10^{-9}$, we can clearly observe an extended inertial range, as it is named in the phenomenology of turbulence, where the spectrum exhibits a power-law behavior whose exponent is governed by the parameter $H$ (here, recall that we chose $H=1/3$), consistently with the prediction obtained for the fractional Gaussian case (Eq. \ref{eq:LimInfTimeExprCorrUXFracCase}). We remark also that our numerical results for the intermittent and non-Gaussian situation ($\gamma=\sqrt{0.02}$) are indistinguishable in this range, as it is expected from the inspection of the observed independence of the series of Figs. \ref{fig:SecondOrder}(a) and (b) to the explored values of $\gamma$.
This is a possible consequence of the fact that the intermittency coefficient has a small value, comparable to those observed in experiments and in simulations of the Navier-Stokes equations \cite{Fri95,CheCas12}. Only in the dissipative range, that is for wavelengths bigger than the characteristic viscous one $k_\nu$ (see Remark \ref{Rem:AddNuFracDyn}), power spectra with different intermittency coefficient $\gamma$ differ. We superimpose with a black line the prediction obtained in the inviscid case ($\nu=0$) which is provided in Eq. \ref{eq:LimInfTimeExprCorrUXFracCase}. Notice that we could have computed in an exact fashion the remaining integral entering in Eq. \ref{eq:LimInfTimeExprCorrUXFracCase} using the expression of the spatial correlation of the force (Eq. \ref{eq:NumDefHatCf}) and special functions, we perform instead a convenient numerical integration. To take into account some implications of the inhomogeneity induced by the truncated version of the force $f_{{\text{trunc}}}$ (Eq. \ref{eq:DefFtrunc}), we propose to weigh the prediction made in Eq. \ref{eq:LimInfTimeExprCorrUXFracCase} by a multiplicative factor given by the integral of the square of the windowing function that enters in its definition. This corresponds to the fraction of energy that is subtracted from the system by the truncation. Accordingly, this factor is defined by and evaluated numerically as $\int \exp[-2x^2/( L^2_{{\text{tot}}}/4-x^2)]dx\approx 0.49 L_{\text{tot}}$, where the integration is made over $|x|< L_{{\text{tot}}}/2$. We observe in the inertial range a nearly perfect collapse of data and prediction, even when $\gamma\ne 0$.

Similarly as for Fig. \ref{fig:SecondOrder}(c), we display in Fig. \ref{fig:SecondOrder}(d) the corresponding second-order structure function $\E|\delta_\ell u_{H,\gamma,\nu}|^2$ as a function of the scale $\ell$, in a logarithmic representation, for the same set of data used in   Fig. \ref{fig:SecondOrder}(c). To estimate this expectation, we average the square norm of the increment $\delta_\ell u_{H,\gamma,\nu}(t,x)=u_{H,\gamma,\nu}(t,x+\ell)- u_{H,\gamma,\nu}(t,x)$ over several independent instances of the solution in time, and also over the region of space $x/L_{{\text{tot}}}\in]-0.2,0.2[$ in which the solution is statistically homogeneous to a good approximation. At large scales, i.e. greater than the integral length scale $L$, the increment variance reaches a plateau, barely dependent on viscosity, which coincides with twice the variance of the solution. Once again, we observe, as $\nu$ decreases, the development of an inertial range where the  second-order structure function behaves as a power-law, whose exponent is governed by the parameter $H$, in a consistent manner with the power-law behavior of the power-spectrum in the corresponding range of wavelengths (Fig. \ref{fig:SecondOrder}(c)), although the power-law is not as clear. Nonetheless, as $\nu$ decreases, we can see this behavior gets closer to the asymptotic prediction that we presented in Eq. \ref{eq:S2FracCaseInviscid} and that we superimpose with a straight black line in Fig. \ref{fig:SecondOrder}(c), weighted for the same reason as for the power spectrum by the factor $0.49 L_{\text{tot}}$. At smaller scales than those pertaining to the inertial range, we recover a scaling behavior proportional to $\ell^2$, as a consequence of the differentiability of the solution for finite viscosity.

Let us make the important remark that, whereas the multifractal ansatz
$v_{H,\gamma}$ (Eq. \ref{eq:VHGammaDef}) exhibits an intermittent correction on the second-order structure function (i.e. take $q=1$ in Eq. \ref{eq:HOStatVHGamma}), it does not seem to be the case from a numerical point of view for the solution of Eq. \ref{eq:PropNumPDE}. This is most certainly related to the fact that we are not presently studying a dynamical version of  $v_{H,\gamma}$ (Eq. \ref{eq:VHGammaDef}), whose evolution is not obviously closed (see the devoted discussion in Paragraph \ref{Par:InducedDynComplexMC}), but an approximate evolution that has been obtained following a closure approach (see Paragraph \ref{Par:ClosureInducedDynComplexMC}).

\begin{figure}[t]
    \centering
 \includegraphics[width=\linewidth]{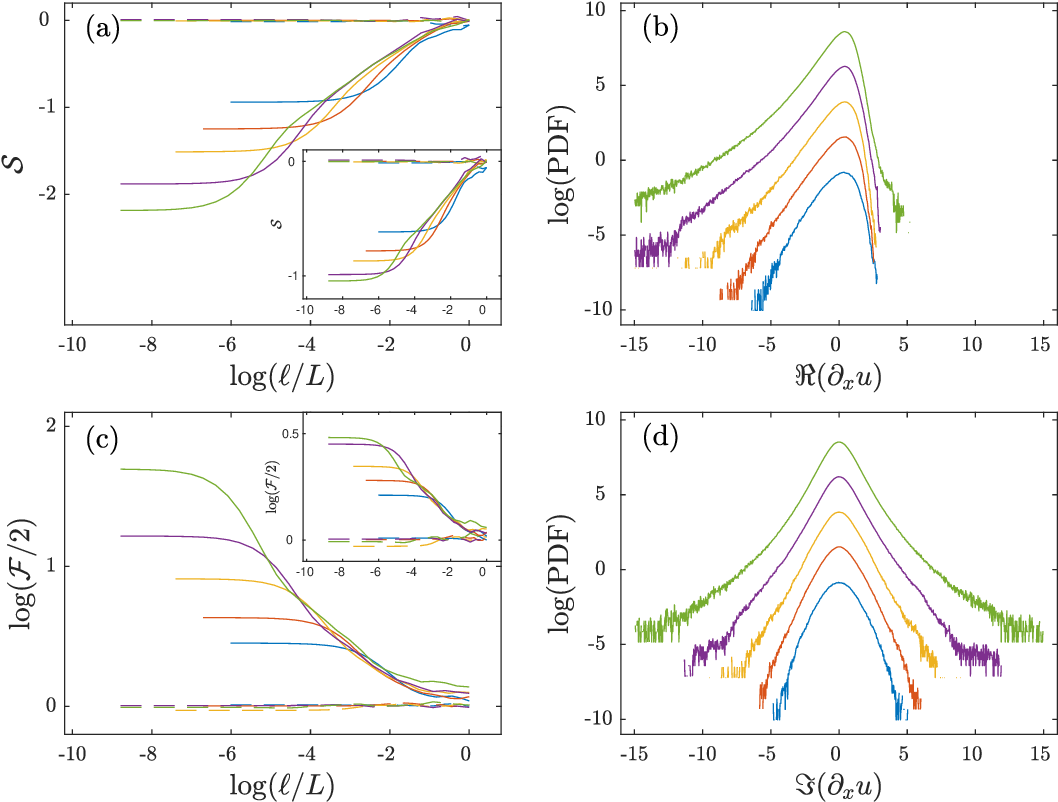}
    \caption{(a) and (c): similar plots as in Fig. \ref{fig:SecondOrder}(d), for the same values of the parameters, using the same colors and representing $\gamma=0$ with dashed lines, and $\gamma=\sqrt{0.02}$ with solid ones, but concerning the increment skewness $\mathcal S(\ell)$ (Eq. \ref{eq:DefSkew}) in (a), and flatness $\mathcal F(\ell)$ (Eq. \ref{eq:DefFlat}) in (c). Corresponding insets show the same statistical quantities but for a smaller value of the parameter $\gamma=\sqrt{0.01}$. (b) and (d): histograms of the values of the real $\Re \partial_xu_{H,\gamma,\nu}$ and imaginary $\Im \partial_xu_{H,\gamma,\nu}$ parts of the derivative of the solution, with the same parameters as former figures but for the single value $\gamma=\sqrt{0.02}$. For the sake of clarity, histograms correspond to unit-variance probability density functions, and are arbitrary vertically shifted (see text).}\label{fig:HOStat}
  \end{figure}

Let us now discuss higher-order statistics than the second-order one, and thus quantify the effects of the quadratic term entering in Eq. \ref{eq:PropNumPDE}, which introduces the parameter $\gamma$ in the picture. Let us first introduce the skewness $\mathcal S(\ell)$ of the increments, i.e.
\begin{align}\label{eq:DefSkew}
\mathcal S(\ell) = \frac{\E\left(\Re\delta_\ell u_{H,\gamma,\nu}\right)^3 }{\left[\E\left(\Re\delta_\ell u_{H,\gamma,\nu}\right)^2\right]^{3/2}},
\end{align}
where only enters the real part of the increment. We display in Fig. \ref{fig:HOStat}(a) the skewness factor of the real part of the increment (Eq. \ref{eq:DefSkew}) as a function of the logarithm of the scales $\ell$, with the same colors representing various values of viscosity as in Fig. \ref{fig:SecondOrder}, using dashed lines for $\gamma=0$ and solid lines for $\gamma=\sqrt{0.02}$. To estimate the expectations entering in Eq. \ref{eq:DefSkew}, we use the same averaging procedure as it is detailed while discussing the results of Fig. \ref{fig:SecondOrder}(d). We indeed observe that $\mathcal S$ vanishes at any scale for $\gamma=0$, consistently with the expected skewness of Gaussian processes. For $\gamma=\sqrt{0.02}$, the scale dependence is rather different. For a given value of viscosity $\nu$, the skewness is negative below the integral scale $\ell\le L$. It then saturates in the dissipative range to converge towards the skewness of the real part of the derivative $\partial_xu_{H,\gamma,\nu}$. As $\nu$ decreases, it seems that the evolution of $\mathcal S$ towards larger negative values follows an approximately viscosity-independent curve in the inertial range. It is rather difficult to see a power law behavior, especially in this representation, but we can say that, while inspecting the numerical results with a smaller value for the parameter $\gamma=\sqrt{0.01}$ as displayed in the inset of Fig. \ref{fig:HOStat}(a), if power-law there is, then the data are compatible with an exponent two times smaller. This suggests that the power law exponent depends quadratically on $\gamma$, in a similar way as in Eq. \ref{eq:HOStatVHGamma}, which was derived for the proposed multifractal ansatz $v_{H,\gamma}$ (Eq. \ref{eq:VHGammaDef}). We could also have computed the skewness (Eq. \ref{eq:DefSkew}) based on the imaginary part of the increments, instead of the real one. In this case, we obtain a vanishing skewness at all scales, even when $\gamma\ne 0$ (data not shown).

We display the scale dependence of the flatness factor of increments in Fig. \ref{fig:HOStat}(c), i.e.
\begin{align}\label{eq:DefFlat}
\mathcal F(\ell) = \frac{\E\left|\delta_\ell u_{H,\gamma,\nu}\right|^4 }{\left[\E\left|\delta_\ell u_{H,\gamma,\nu}\right|^2\right]^{2}},
\end{align}
and we use the same colors for various viscosities and dashed and solid lines respectively for $\gamma=0$ and $\gamma=\sqrt{0.02}$, as we did in Figs. \ref{fig:SecondOrder} and \ref{fig:HOStat}(c). We first observe that $\mathcal F(\ell)\approx 2$ at any scale $\ell$ when $\gamma=0$, as is expected from a complex Gaussian random field, whose real and imaginary parts are independent. For the more interesting case $\gamma=\sqrt{0.02}$, we observe that the flatness departs from the Gaussian value $2$ as $\ell\le L$. In the inertial range of scales, $\mathcal F$ seems to behave as a power law, independently of the value of viscosity. This is a characteristic feature of multifractal processes, in particular reproduced by our multifractal ansatz $v_{H,\gamma}$ (Eq. \ref{eq:VHGammaDef}). Similarly to the skewness, the power law exponent of this observed behavior is tricky to understand. By inspection of the behavior of flatnesses for a smaller intermittency parameter $\gamma=\sqrt{0.01}$ displayed in the inset of Fig. \ref{fig:HOStat}(c), we can infer that again data are compatible with a power law exponent proportional to $\gamma^2$, as is expected from multifractal processes (Eq. \ref{eq:HOStatVHGamma}). The multiplicative factor in front of $\gamma^2$ remains difficult to determine at this stage, due to the limited inertial range observed. A further numerical investigation on the effects of different spatial forcing correlations and of the closure of Eq.~\ref{eq:ApproxW} on these nonlinear exponents is a future perspective, in order to verify the universality properties predicted by Eqs.~\ref{eq:HOStatVHGamma} and \ref{eq:OddStatVHGamma}.

Finally, we display respectively in Fig. \ref{fig:HOStat}(b) and (d) the histograms of the real and imaginary parts of the gradients $\partial_xu_{H,\gamma,\nu}$ for various values of viscosity and $\gamma=\sqrt{0.02}$, using the same colors as those that have been used formerly. This estimation of the probability density functions (PDFs) is made following the same averaging procedure, that is over independent instances in time and across the approximately statistically homogeneous region $x/L_{{\text{tot}}}\in]-0.2,0.2[$. To make the comparison clear between different viscosities, we display the estimated PDFs such that they are all of unit variance, and we shift them vertically in an arbitrary manner to highlight the evolution of their shape. As expected, for $\gamma=0$, PDFs of gradients are Gaussian for any viscosity (data not shown). On the contrary, for $\gamma=\sqrt{0.02}$, we observe a continuous deformation of their shape as $\nu$ decreases, being closer to a Gaussian shape at high viscosity $\nu=10^{-5}$, and exhibiting wider and wider tails as $\nu$ decreases towards its lowest value $\nu=10^{-9}$. Consistently with the observed behavior of the skewness $\mathcal S$ (Eq. \ref{eq:DefSkew}), PDFs are negatively skewed for $\Re\partial_xu_{H,\gamma,\nu}$ and symmetrical for $\Im\partial_xu_{H,\gamma,\nu}$, as it is obtained for the multifractal ansatz $v_{H,\gamma}$ (Eq. \ref{eq:VHGammaDef}), and in particular in its third-order structure function (Eq.~\ref{eq:OddStatVHGamma}). Also, consistently with the fact that the small scale plateaus in the increment flatnesses rise as the viscosity decreases (Fig. \ref{fig:HOStat}(c)), wider tail of the estimated gradients PDFs develop for smaller viscosities.

\appendix

\section{Proof of Proposition \ref{Def:MultiFieldV}}\label{App:ProofVHGamma}

Let us start with proposing Lemma \ref{Lemma:GIPComplex}, similarly to the Lemma 2.1 of Ref. \cite{RobVar08}, but here for complex Gaussian variables:
\begin{lemma} \label{Lemma:GIPComplex}
 Consider a complex zero average Gaussian random variable $Z$, a function $F:\C\rightarrow \C$ and its derivative $F'$ that grows at most exponentially. We have
 \begin{align}
 \E\left[ZF(Z)\right]=\E(Z^2)\E\left[F'(Z)\right].
 \end{align}
More generally, considering the collection of $(n+1)$ complex Gaussian variables $(Z,Z_1,...,Z_n)$ and the function $F:\C^n\rightarrow \C$. We have the following Gaussian integration by parts formula
 \begin{align}\label{eq:IntByPartsFormula}
 \E\left[ZF(Z_1,...,Z_n)\right]=\sum_{k=1}^n\E(ZZ_k)\E\left[\frac{\partial F}{\partial x_k}(Z_1,...,Z_n)\right].
 \end{align}
  \end{lemma}

 \textit{Proof of Proposition \ref{Def:MultiFieldV}:}
 Concerning the average of $v_{H,\gamma}$ (Eq. \ref{eq:VHGammaDef}), we make use of Eq. \ref{eq:IntByPartsFormula} and obtain
  \begin{align}
\E\left[e^{\gamma \tilde{P}_0u_0(t,y)} u_0(t,y)\right]=\gamma \E\left[ u_0(t,y)\tilde{P}_0u_0(t,y)\right]e^{\frac{\gamma^2}{2}\E\left[\tilde{v}_0^2(t,y)\right]}=0,
  \end{align}
because for any positions, $\E\left[ u_0(t,y)u_0(t,z)\right]=0$ (Eq. \ref{eq:ExprCorrUReIm}), which shows that $\E v_{H,\gamma}=0$.

 The calculation of the variance is done in a similar way, and requires the following step: Making use of Eq. \ref{eq:IntByPartsFormula}, we obtain
   \begin{align}\label{eq:CompCorrOrder2}
&\E\left[ u_0(t,y_1) u_0^*(t,y_2)e^{\gamma \left(\tilde{P}_0u_0(t,y_1)+\tilde{P}_0u_0^*(t,y_2)\right)}\right]=\mathcal C_{u_0}(t,y_1-y_2)e^{\gamma^2\mathcal C_{\tilde{v}_0}(t,y_1-y_2)}\\
&+\gamma \E\left[ u_0(t,y_1) \tilde{P}_0u_0^*(t,y_2)\right]\E\left[ u_0^*(t,y_2)e^{\gamma \left(\tilde{P}_0u_0(t,y_1)+\tilde{P}_0u_0^*(t,y_2)\right)}\right]\notag\\
&= \left(\mathcal C_{u_0}(t,y_1-y_2) + \gamma^2 \E\left[ u_0(t,y_1) \tilde{P}_0u_0^*(t,y_2)\right] \E\left[ u_0^*(t,y_2) \tilde{P}_0u_0(t,y_1)\right]\right)e^{\gamma^2\mathcal C_{\tilde{v}_0}(t,y_1-y_2)}.\notag
  \end{align}
Using the odd symmetry of the function $ \tilde{P}_0(x)=- \tilde{P}_0(-x)$, notice that
   \begin{align}\label{eq:DefKOdd}
\mathcal K(t,y_1-y_2)\equiv\E\left[ u_0(t,y_1) \tilde{P}_0u_0^*(t,y_2)\right]  &=\int  \tilde{P}_0(y_2-z)\mathcal C_{u_0}(t,y_1-z)dz\\
&=-\E\left[ u_0^*(t,y_2) \tilde{P}_0u_0(t,y_1)\right] \notag\\
&=-\mathcal K^*(t,y_2-y_1) \notag
  \end{align}
and we obtain
   \begin{align}\label{eq:VarFiniteTimeVH}
\E\left[\left|v_{H,\gamma}(t,x)\right|^2\right]=  \int e^{2i\pi ky}\frac{1}{|k|_{1/L}^{2H+1}}\left[ \mathcal C_{u_0}(t,y)-\gamma^2 \mathcal K^2(t,y)\right]e^{\gamma^2\mathcal C_{\tilde{v}_0}(t,y)}dkdy.
   \end{align}
Remark now that
   \begin{align}\label{eq:VarFiniteTimeVHRemark}
   \int_0^{+\infty}(P_H\star P_H)'(y)e^{\gamma^2\mathcal C_{\tilde{v}_0}(t,y)}dy = &-(P_H\star P_H)(0)e^{\gamma^2\mathcal C_{\tilde{v}_0}(t,0)} \\
 &  -   \int_0^{+\infty}(P_H\star P_H)(y)\gamma^2\mathcal C'_{\tilde{v}_0}(t,y)e^{\gamma^2\mathcal C_{\tilde{v}_0}(t,y)}dy,\notag
   \end{align}
where we have introduced the correlation product $\star$, defined by, for any appropriate real functions $f$ and $g$,
   \begin{align}\label{eq:DefCorrProductStar}
(f\star g)(y)=\int f(x)g(x+y)dx =   \int e^{2i\pi ky} \hat{f}^*(k)\hat{g}(k)dk,
   \end{align}
such that
   \begin{align}
&\E\left[\left|v_{H,\gamma}(t,x)\right|^2\right] + \frac{\mathcal C_f(0)}{2|c|}\int_0^{+\infty}(P_H\star P_H)'(y)e^{\gamma^2\mathcal C_{\tilde{v}_0}(t,y)}dy=\label{eq:PrincContribVar}\\
&\int (P_H\star P_H)(y) \mathcal C_{u_0}(t,y) e^{\gamma^2\mathcal C_{\tilde{v}_0}(t,y)}dy -\frac{\mathcal C_f(0)}{2|c|}(P_H\star P_H)(0)e^{\gamma^2\mathcal C_{\tilde{v}_0}(t,0)}\label{eq:FirstContribVar}\\
&-\gamma^2 \int (P_H\star P_H)(y) \left[ \mathcal K^2(t,y)+\frac{\mathcal C_f(0)}{4|c|}\mathcal C'_{\tilde{v}_0}(t,|y|)\right]e^{\gamma^2\mathcal C_{\tilde{v}_0}(t,y)}dy.\label{eq:SecondContribVar}
   \end{align}
Recall that $\mathcal C_{u_0}(t,y) $ behaves as a Dirac function as $t\to\infty$, weighted by an appropriate factor, as is stated in Eq. \ref{eq:BehavCuWN}. Hence, it is clear that the contribution given in Eq. \ref{eq:FirstContribVar} will vanish as $t\to\infty$. Similarly, in the same limit, the function $\mathcal K(t,y)$ (Eq. \ref{eq:DefKOdd}) converges towards $\tilde{P}_0(-y)$, again weighted by an appropriate factor, such that, using the asymptotic expression of $\mathcal C_{\tilde{v}_0}$ (Eq. \ref{eq:CovarVFGFH0DefHTilde}), we obtain pointwise
   \begin{align}\label{eq:AsymptKVar}
\lim_{t\to\infty} \left[\mathcal K^2(t,y)+\frac{\mathcal C_f(0)}{4|c|}\mathcal C'_{\tilde{v}_0}(t,|y|)\right]=\frac{\mathcal C^2_f(0)}{4|c|^2}\left[\tilde{P}_0^2(y)-|\ln_+'(|y|)| \right]+\frac{\mathcal C_f(0)}{4|c|}\tilde{h}'(|y|),
   \end{align}
where we have denoted by $|\ln_+'(|y|)| $ the derivative of the smoothly-truncaded logarithm $\ln_+$ evaluated at $|y|$, that is expected to behave as $1/|y|$ in the vicinity of the origin.

Let us focus on the second contribution displayed in Eq. \ref{eq:SecondContribVar}. As we have already observed, the function $(P_H\star P_H)(y)$ is a bounded function of its argument for any $H>0$, and its rapid decrease away from the origin ensures integrability when the dummy variable $y$ goes towards infinite values. Notice that
\begin{align}
\tilde{P}_0(y) =-i \int e^{2i\pi ky}\frac{k}{|k|^{3/2}_{1/L}}dk\build{\sim}_{y\to 0}^{}\frac{y}{|y|^{3/2}},
\end{align}
which says that the first term of the RHS of Eq. \ref{eq:AsymptKVar} grows at most logarithmically near the origin, which is integrable. Thus, the integral entering in Eq. \ref{eq:SecondContribVar} exists as $t\to\infty$ if the remaining singular term, i.e. $|y|^{-\gamma^2\mathcal C_f(0)/|c|}$, is integrable, i.e. $\gamma^2<|c|/\mathcal C_f(0)$.

To conclude, concerning the limit at large time of the variance $\E\left[\left|v_{H,\gamma}(t,x)\right|^2\right]$, let us examine the second term of the LHS of Eq. \ref{eq:PrincContribVar}. It is easy to see that near the origin, whereas $(P_H\star P_H)(y)$ remains bounded for any $H>0$, its derivative will behave as $(P_H\star P_H)'(y)\sim y^{2H-1}$. Thus, as $t\to \infty$, this contribution is finite for $2H-1-\gamma^2\mathcal C_f(0)/|c|>-1$, i.e. $\gamma^2<2H|c|/\mathcal C_f(0)$. Hence, for $\gamma^2<\min(1,2H)|c|/\mathcal C_f(0)$, the variance is finite and its expression is given by
   \begin{align}\label{eq:LimVarVHgamma}
&\lim_{t\to\infty}\E\left[\left|v_{H,\gamma}(t,x)\right|^2\right] =- \frac{\mathcal C_f(0)}{2|c|}\int_0^{+\infty}(P_H\star P_H)'(y)\left|\frac{L}{y}\right|_+^{\gamma^2\frac{\mathcal C_f(0)}{|c|}}e^{\gamma^2\tilde{h}(y)}dy\\
&-\gamma^2\frac{\mathcal C_f(0)}{4|c|} \int (P_H\star P_H)(y) \left[ \frac{\mathcal C_f(0)}{|c|}\left[\tilde{P}_0^2(y)-|\ln_+'(|y|)| \right]+\tilde{h}'(|y|)\right]\left|\frac{L}{y}\right|_+^{\gamma^2\frac{\mathcal C_f(0)}{|c|}}e^{\gamma^2\tilde{h}(y)}dy,\notag
\end{align}
where the notation $|y|_+=\exp\left( \ln_+|y|\right)$ is introduced.

To see the behavior at small scales of the second-order structure function, consider the function
\begin{align}\label{eq:DefMatCalPHIncr}
\mathcal P_{H,\ell}(x)=P_{H}(x+\ell)-P_{H}(x)=\int e^{2i\pi kx}\frac{e^{2i\pi k\ell}-1}{|k|_{1/L}^{H+1/2}}dk,
\end{align}
such that we can conveniently write the velocity increment as
\begin{align}\label{eq:DefIncrVHGamma}
\delta_\ell v_{H,\gamma}(t,x)=\int \mathcal P_{H,\ell}(x-y)e^{\gamma \tilde{P}_0u_0(t,y)} u_0(t,y)dy.
\end{align}
Previous calculations concerning the variance apply and we get
   \begin{align}\label{eq:LimVarDeltaLVHgamma}
&\lim_{t\to\infty}\E\left[\left|\delta_\ell v_{H,\gamma}(t,x)\right|^2\right] =- \frac{\mathcal C_f(0)}{2|c|}\int_0^{+\infty}( \mathcal P_{H,\ell}\star  \mathcal P_{H,\ell})'(y)\left|\frac{L}{y}\right|_+^{\gamma^2\frac{\mathcal C_f(0)}{|c|}}e^{\gamma^2\tilde{h}(y)}dy\\
&-\gamma^2\frac{\mathcal C_f(0)}{4|c|} \int ( \mathcal P_{H,\ell}\star  \mathcal P_{H,\ell})(y) \left[ \frac{\mathcal C_f(0)}{|c|}\left[\tilde{P}_0^2(y)-|\ln_+'(|y|)| \right]+\tilde{h}'(|y|)\right]\left|\frac{L}{y}\right|_+^{\gamma^2\frac{\mathcal C_f(0)}{|c|}}e^{\gamma^2\tilde{h}(y)}dy.\notag
\end{align}
Notice that
\begin{align}\label{eq:DefConvolMatCalPHIncr}
( \mathcal P_{H,\ell}\star  \mathcal P_{H,\ell})(y)=\int e^{2i\pi ky}\frac{\left|e^{2i\pi k\ell}-1\right|^2}{|k|_{1/L}^{2H+1}}dk,
\end{align}
such that
\begin{align}\label{eq:RescaleConvolMatCalPHIncr}
( \mathcal P_{H,\ell}\star  \mathcal P_{H,\ell})(\ell y)\build{\sim}_{\ell\to 0^+}^{}\ell^{2H}\int e^{2i\pi ky}\frac{\left|e^{2i\pi k}-1\right|^2}{|k|^{2H+1}}dk,
\end{align}
this equivalence at small scales making sense only for $H\in]0, 1[$. Similarly, we have
\begin{align}\label{eq:RescaleDerivConvolMatCalPHIncr}
( \mathcal P_{H,\ell}\star  \mathcal P_{H,\ell})'(\ell y)\build{\sim}_{\ell\to 0^+}^{}\ell^{2H-1}g_H(y),
\end{align}
where
\begin{align}\label{eq:DefGHRescaleDerivConvolMatCalPHIncr}
g_H(y)=\int 2i\pi k e^{2i\pi ky}\frac{\left|e^{2i\pi k}-1\right|^2}{|k|^{2H+1}}dk.
\end{align}
Once having rescaled the dummy variable $y$ entering in the integrals at the RHS of Eq. \ref{eq:LimVarDeltaLVHgamma}, we can see that the first term will be order $\ell^{2H-\gamma^2\mathcal C_f(0)/|c|}$, and thus will dominate the second term that goes to zero as $\ell^{2H+1-\gamma^2\mathcal C_f(0)/|c|}$. Doing so, we get the equivalent behavior of the second-order structure function at small scales, which reads
 \begin{align}\label{eq:EquivSmallLLimVarDeltaLVHgamma}
\lim_{t\to\infty}\E\left[\left|\delta_\ell v_{H,\gamma}(t,x)\right|^2\right] \build{\sim}_{\ell\to 0^+}^{}- \frac{\mathcal C_f(0)}{2|c|}\ell^{2H}\left(\frac{\ell}{L}\right)^{-\gamma^2\frac{\mathcal C_f(0)}{|c|}}e^{\gamma^2\tilde{h}(0)}\int_0^{+\infty}\frac{g_H(y)}{|y|^{\gamma^2\frac{\mathcal C_f(0)}{|c|}}}dy.
\end{align}
It remains to determine the range of parameters such that the equivalence given in Eq. \ref{eq:EquivSmallLLimVarDeltaLVHgamma} makes sense, and hence check the integrability of the remaining integral that enters in it. Although the behavior of the function $g_H$ defined in Eq. \ref{eq:DefGHRescaleDerivConvolMatCalPHIncr} at small and large arguments can be tricky to establish, its integrability is pretty much straightforward. Indeed, with notation $a=\gamma^2\frac{\mathcal C_f(0)}{|c|}$, using the equality
$$ \int_0^{+\infty}\frac{g_H(y)}{|y|^{a}}dy=- (2\pi)^{a}\Gamma\left(1 - a\right) \cos\left(\frac{a\pi}{2}\right)\int \frac{\left|e^{2i\pi k}-1\right|^2}{|k|^{2H+1-a}}dk,$$
it is clear that the equivalent (Eq. \ref{eq:EquivSmallLLimVarDeltaLVHgamma}) makes sense for $H\in]0,1[$ and $\gamma^2\frac{\mathcal C_f(0)}{|c|}<\min(1,2H)$, and is indeed positive. As a further check, we can note that the expression  (Eq. \ref{eq:EquivSmallLLimVarDeltaLVHgamma}) indeed coincides with the equivalence obtained for fractional Gaussian fields (Eq. \ref{eq:CalcVarIncrFGF}) when $\gamma=0$.

Let us now calculate the third order structure function. We have, making use of the definition and symmetries of the function $\mathcal K$ (Eq. \ref{eq:DefKOdd}),
\begin{align}\label{eq:CompCorrOrder3}
&\E\left[ u_0(t,y_1) u_0^*(t,y_2)u_0(t,y_3)e^{\gamma \left(\tilde{P}_0u_0(t,y_1)+\tilde{P}_0u_0^*(t,y_2)+\tilde{P}_0u_0(t,y_3)\right)}\right]\\
&=\mathcal C_{u_0}(t,y_1-y_2)\E\left[ u_0(t,y_3)e^{\gamma \left(\tilde{P}_0u_0(t,y_1)+\tilde{P}_0u_0^*(t,y_2)+\tilde{P}_0u_0(t,y_3)\right)}\right]\notag\\
&+\gamma\mathcal K(t,y_1-y_2)\E\left[ u_0^*(t,y_2)u_0(t,y_3)e^{\gamma \left(\tilde{P}_0u_0(t,y_1)+\tilde{P}_0u_0^*(t,y_2)+\tilde{P}_0u_0(t,y_3)\right)}\right]\notag\\
&=\mathcal C_{u_0}(t,y_1-y_2)\left[ \gamma\mathcal K(t,y_3-y_2)e^{\gamma^2\left[\mathcal C_{\tilde{v}_0}(t,y_1-y_2)+\mathcal C_{\tilde{v}_0}(t,y_3-y_2)\right] }\right]\notag\\
&+\gamma \mathcal K(t,y_1-y_2)\left[ \mathcal C_{u_0}(t,y_3-y_2)+\gamma^2\left[\mathcal K^*(t,y_2-y_1)+\mathcal K^*(t,y_2-y_3)\right]\mathcal K(t,y_3-y_2)\right]\notag\\
&\times e^{\gamma^2\left[\mathcal C_{\tilde{v}_0}(t,y_1-y_2)+\mathcal C_{\tilde{v}_0}(t,y_3-y_2)\right] } \notag\\
&=\gamma \left[\mathcal C_{u_0}(t,y_1-y_2)\mathcal K(t,y_3-y_2)+\mathcal C_{u_0}(t,y_3-y_2)\mathcal K(t,y_1-y_2)\right]e^{\gamma^2\left[\mathcal C_{\tilde{v}_0}(t,y_1-y_2)+\mathcal C_{\tilde{v}_0}(t,y_3-y_2)\right] }\notag\\
&-\gamma^3\mathcal K(t,y_1-y_2)\left[\mathcal K(t,y_1-y_2)+\mathcal K(t,y_3-y_2)\right]\mathcal K(t,y_3-y_2)   e^{\gamma^2\left[\mathcal C_{\tilde{v}_0}(t,y_1-y_2)+\mathcal C_{\tilde{v}_0}(t,y_3-y_2)\right] },\notag
  \end{align}
such that
\begin{align*}
\E \left[ \delta_\ell v_{H,\gamma}|\delta_\ell v_{H,\gamma}|^2\right]&=2\gamma \int \mathcal P_{H,\ell}(y_2-z_1) \mathcal P_{H,\ell}(y_2) \mathcal P_{H,\ell}(y_2-z_3) \\
&\times\left[ \mathcal C_{u_0}(t,z_1)-\gamma^2\mathcal K^2(t,z_1)\right]\mathcal K(t,z_3)e^{\gamma^2\left[\mathcal C_{\tilde{v}_0}(t,z_1)+\mathcal C_{\tilde{v}_0}(t,z_3)\right] }dz_1 dy_2 dz_3.\notag
\end{align*}
Let us introduce the following function
 \begin{align*}
h_{H,\ell}(z_1,z_3)&=\int\mathcal P_{H,\ell}(y_2-z_1) \mathcal P_{H,\ell}(y_2) \mathcal P_{H,\ell}(y_2-z_3) dy_2\\
&=\int  e^{-2i\pi (k_1z_1+k_3z_3)} \frac{\left(e^{2i\pi k_1\ell}-1\right)\left(e^{-2i\pi (k_1+k_3)\ell}-1\right)\left(e^{2i\pi k_3\ell}-1\right)}{|k_1|_{1/L}^{H+1/2}|k_1+k_3|_{1/L}^{H+1/2}|k_3|_{1/L}^{H+1/2}}dk_1dk_3\\
&=-h_{H,\ell}(-z_1,-z_3),
\end{align*}
such that
\begin{align}\label{eq:DefMom3ProcessFiniteTime}
&\E \left[ \delta_\ell v_{H,\gamma}|\delta_\ell v_{H,\gamma}|^2\right]=\\
&2\gamma \int h_{H,\ell}(z_1,z_3) \left[ \mathcal C_{u_0}(t,z_1)-\gamma^2\mathcal K^2(t,z_1)\right]\mathcal K(t,z_3)e^{\gamma^2\left[\mathcal C_{\tilde{v}_0}(t,z_1)+\mathcal C_{\tilde{v}_0}(t,z_3)\right] }dz_1dz_3\notag.
\end{align}
Using the same ideas to determine the limiting value as $t\to\infty$ of the variance (Eq. \ref{eq:VarFiniteTimeVH}), remark that
\begin{align}\label{eq:VarFiniteTimeVHRemarkThird}
   \int_0^{+\infty} \partial_{z_1}h_{H,\ell}(z_1,z_3)e^{\gamma^2\mathcal C_{\tilde{v}_0}(t,z_1)}dz_1 = &-h_{H,\ell}(0,z_3)e^{\gamma^2\mathcal C_{\tilde{v}_0}(t,0)} \\
 &  -   \int_0^{+\infty}h_{H,\ell}(z_1,z_3)\gamma^2\mathcal C'_{\tilde{v}_0}(t,z_1)e^{\gamma^2\mathcal C_{\tilde{v}_0}(t,z_1)}dz_1,\notag
 \end{align}
as we obtained in Eq. \ref{eq:VarFiniteTimeVHRemark}. Doing so, we determine the proper quantity that eventually dominates at small scales, and we obtain
\begin{align*}
\lim_{t\to\infty}\E \left[ \delta_\ell v_{H,\gamma}|\delta_\ell v_{H,\gamma}|^2\right]&\build{\sim}_{\ell\to 0^+}^{}\gamma \frac{\mathcal C^2_f(0)}{2|c|^2}\ell^{3H}\left(\frac{\ell}{L}\right)^{-2\gamma^2\frac{\mathcal C_f(0)}{|c|}}e^{2\gamma^2\tilde{h}(0)} \\
&\times\int_{(z_1,z_3)\in \R^+\times\R}g_H(z_1,z_3)\frac{1}{|z_1|^{\frac{\gamma^2\mathcal C_f(0)}{|c|}}}\frac{z_3}{|z_3|^{\frac{3}{2}+\frac{\gamma^2\mathcal C_f(0)}{|c|}}}dz_1dz_3,
\end{align*}
where we have introduced the function
\begin{align}\label{eq:DefGHZ1Z2}
g_H(z_1,z_3)=-2i\pi\int  e^{-2i\pi (k_1z_1+k_3z_3)} \frac{k_1\left(e^{2i\pi k_1}-1\right)\left(e^{-2i\pi (k_1+k_3)}-1\right)\left(e^{2i\pi k_3}-1\right)}{|k_1|^{H+1/2}|k_1+k_3|^{H+1/2}|k_3|^{H+1/2}}dk_1dk_3.
\end{align}
Additionally, we will need the following exact Fourier transforms,
\begin{align}\label{eq:FT1Comp}
 \int  e^{-2i\pi k_1z_1} \frac{1_{z_1\ge 0}}{|z_1|^a}dz_1 =(2\pi)^{a-1}  \Gamma(1 - a) \left[\sin(a \pi/2) -i\cos(a \pi/2)\frac{k_1}{|k_1|}\right]\frac{1}{|k_1|^{1-a}},
 \end{align}
for $0<a<1$, and
\begin{align}\label{eq:FT2Comp}
\int  e^{-2i\pi k_3z_3} \frac{z_3}{|z_3|^{\frac{3}{2}+a}}dz_3 =-i (2 \pi)^{a+1/2}
   \frac{1/4 + a/2 }{\Gamma(a+3/2)\sin(\pi (a/2+1/4))} \frac{k_3}{|k_3|^{3/2 - a}},
    \end{align}
for $0<a<3/2$, and the identity
\begin{align*}\left(e^{2i\pi k_1}-1\right)\left(e^{-2i\pi (k_1+k_3)}-1\right)&\left(e^{2i\pi k_3}-1\right)\\
&=-2i\left[ \sin\left(2\pi(k_1+k_3)\right)-\sin(2\pi k_1)-\sin(2\pi k_3)\right].
\end{align*}

Using symmetries, it can be shown that the real part of Eq. \ref{eq:FT1Comp} does not contribute, only remaining
\begin{align}\label{eq:IntGZ1Z3AgainstMC}
\int_{(z_1,z_3)\in \R^+\times\R}g_H(z_1,z_3)&\frac{1}{|z_1|^{\frac{\gamma^2\mathcal C_f(0)}{|c|}}}\frac{z_3}{|z_3|^{\frac{3}{2}+\frac{\gamma^2\mathcal C_f(0)}{|c|}}}dz_1dz_3\\
&=4\pi A_\gamma \int \frac{k_3\left[ \sin\left(2\pi(k_1+k_3)\right)-\sin(2\pi k_1)-\sin(2\pi k_3)\right]}{|k_1|^{H+1/2-\frac{\gamma^2\mathcal C_f(0)}{|c|}}|k_1+k_3|^{H+1/2}|k_3|^{H+2-\frac{\gamma^2\mathcal C_f(0)}{|c|}}}dk_1dk_3\notag,
\end{align}
with $A_\gamma\in\R$ a real multiplicative constant that can be obtained from the multiplicative contributions displayed in Eqs. \ref{eq:FT1Comp} and  \ref{eq:FT2Comp}. The sign of the remaining contribution of the RHS of Eq. \ref{eq:IntGZ1Z3AgainstMC} expressed as a double integral over the dummy variables $k_1$ and $k_3$ is not obvious, neither whether it vanishes or not. Nonetheless, it gives a condition on $\gamma$, to ensure its finiteness. Inspecting the integrability properties of this term, we find that the integral exists along the diagonal $k_1=k_3$ if
\begin{align}\label{eq:CondM3}
\frac{\gamma^2\mathcal C_f(0)}{|c|}<\min(1,3H/2).
\end{align}
Doing so, we have thus shown that, under the condition provided in Eq. \ref{eq:CondM3}, the third-order moment of the increments of the process $v_{H,\gamma}$, as it is defined in Eq. \ref{eq:DefMom3ProcessFiniteTime}, is finite, and does not vanish in an obvious manner. It is furthermore real, and it behaves at small scale as
\begin{align*}
\lim_{t\to\infty}\E \left[ \delta_\ell v_{H,\gamma}|\delta_\ell v_{H,\gamma}|^2\right]&\build{\sim}_{\ell\to 0^+}^{}d_{H,\gamma}\ell^{3H}\left(\frac{\ell}{L}\right)^{-2\gamma^2\frac{\mathcal C_f(0)}{|c|}}\end{align*}
with
\begin{align}\label{eq:CalcDHGamma3rd}
d_{H,\gamma}=\gamma \frac{\mathcal C^2_f(0)}{2|c|^2}e^{2\gamma^2\tilde{h}(0)} \int_{(z_1,z_3)\in \R^+\times\R}g_H(z_1,z_3)\frac{1}{|z_1|^{\frac{\gamma^2\mathcal C_f(0)}{|c|}}}\frac{z_3}{|z_3|^{\frac{3}{2}+\frac{\gamma^2\mathcal C_f(0)}{|c|}}}dz_1dz_3,
\end{align}
where the function $g_H(z_1,z_3)$ is defined in Eq. \ref{eq:DefGHZ1Z2}.

Let us finally determine the behavior at small scales of the statistics at high-order considering $q\in\N^*$,
\begin{align}\label{eq:Mom2QStart}
\E \left[ |\delta_\ell v_{H,\gamma}|^{2q}\right]&=\int
\prod_{i=1}^q \mathcal P_{H,\ell}(x-y_i)\mathcal P_{H,\ell}(x-z_i)\prod_{i=1}^q dy_idz_i\\
&\times\E\left[\prod_{i=1}^q u_0(t,y_i) u_0^*(t,z_i)e^{\gamma \sum_{i=1}^q\tilde{P}_0u_0(t,y_i)+\tilde{P}_0u_0^*(t,z_i)}\right],\notag
\end{align}
where the operator $\mathcal P_{H,\ell}$ is defined in Eq. \ref{eq:DefMatCalPHIncr}. The determination of the exact expression of the correlator entering in Eq. \ref{eq:Mom2QStart} can be done using some combinatorial analysis, although it can become cumbersome. Instead, in a first approach, let us evaluate the spectrum of exponents that governs the decrease towards 0 as $\ell\to 0$. In particular, intermittent corrections are eventually governed by a term of the form
$$ \E\left[e^{\gamma \sum_{i=1}^q\tilde{P}_0u_0(t,y_i)+\tilde{P}_0u_0^*(t,z_i)}\right] = e^{\gamma^2\left[ \sum_{i=1}^q\mathcal C_{\tilde{v}_0}(t,y_i-z_i)+ \sum_{i<j=1}^q\mathcal C_{\tilde{v}_0}(t,y_i-z_j)+\mathcal C^*_{\tilde{v}_0}(t,y_j-z_i)\right]},$$
contributing at small scales as
\begin{align*}
\lim_{t\to\infty} \E\left[e^{\gamma \sum_{i=1}^q\tilde{P}_0u_0(t,\ell y_i)+\tilde{P}_0u_0^*(t,\ell z_i)}\right] &\build{\sim}_{\ell\to 0^+}^{}\left(\frac{\ell}{L}\right)^{-q^2\gamma^2\frac{\mathcal C_f(0)}{|c|}}e^{q^2\gamma^2\tilde{h}(0)}  \\
&\times \prod_{i=1}^q \frac{1}{|y_i-z_i|^{\gamma^2\frac{\mathcal C_f(0)}{|c|}}} \prod_{i<j=1}^q \frac{1}{|(y_i-z_j)(y_j-z_i)|^{\gamma^2\frac{\mathcal C_f(0)}{|c|}}},
 \end{align*}
whereas contributions from the fractional part will be of the order of $\ell^{2qH}$. Once again, the determination of the appropriate range of values for $\gamma$ is tricky to get at this stage because we have to compute in an exact fashion the expectation entering in the RHS of Eq. \ref{eq:Mom2QStart}. To do so, we have to generalize the calculations made in Eqs. \ref{eq:CompCorrOrder2} and \ref{eq:CompCorrOrder3}, using combinatorial developments such as those proposed in Ref. \cite{RobVar08} (see their Lemma 2.2). Such a calculation is beyond the scope of the present article. We nonetheless expect the additional condition $\gamma^2\frac{\mathcal C_f(0)}{|c|}<2H/q$.

\smallskip
\noindent \textbf{Acknowledgements.} We warmly thank Laure Saint-Raymond, with whom this project has started, for many enlightening discussions. Also, we thank J\'er\'emie Bouttier, Gr\'egory Miermont, R\'emi Rhodes and Simon Thalabard for several discussions on this subject, and Martin Hairer for additional fruitful discussions and for bringing to our knowledge the important Ref. \cite{MatSui07}.  G. B. A. and L. C. are partially supported by the Simons Foundation Award ID: 651475. J.-C. M. is partially supported by the NSF grant DMS-1954357.


%
%



\end{document}